\documentclass[reqno]{amsart}

\usepackage{xypic}
\input xy
\xyoption{all}
\usepackage{epsfig}
\usepackage{color}
\usepackage{amsthm}
\usepackage{amssymb}
\usepackage{amsmath}
\usepackage{amscd}
\usepackage{amsopn}
\usepackage{pinlabel}

\usepackage{url}
\usepackage{hyperref} \hypersetup{colorlinks}

\usepackage{color} 

\definecolor{darkred}{rgb}{1,0,0} 
\definecolor{darkgreen}{rgb}{0,0.8,0}
\definecolor{darkblue}{rgb}{0,0,1}

\hypersetup{colorlinks,
linkcolor=darkblue,
filecolor=darkgreen,
urlcolor=darkred,
citecolor=darkgreen}

\newcommand{\labell}[1] {\label{#1}}

\numberwithin{equation}{section}
\newtheorem {Theorem}{Theorem}
\numberwithin{Theorem}{section}

\newtheorem {Lemma}[Theorem]    {Lemma}

\newtheorem {Proposition}[Theorem]{Proposition}
\newtheorem {Corollary}[Theorem]{Corollary}
\theoremstyle{definition}
\newtheorem{Definition}[Theorem]{Definition}
\theoremstyle{remark}
\newtheorem{Remark}[Theorem]{Remark}
\newtheorem{Example}[Theorem]{Example}

\expandafter\chardef\csname pre amssym.def at\endcsname=\the\catcode`\@
\catcode`\@=11
\def\undefine#1{\let#1\undefined}
\def\newsymbol#1#2#3#4#5{\let\next@\relax
 \ifnum#2=\@ne\let\next@\msafam@\else
 \ifnum#2=\tw@\let\next@\msbfam@\fi\fi
 \mathchardef#1="#3\next@#4#5}
\def\mathhexbox@#1#2#3{\relax
 \ifmmode\mathpalette{}{\m@th\mathchar"#1#2#3}%
 \else\leavevmode\hbox{$\m@th\mathchar"#1#2#3$}\fi}
\def\hexnumber@#1{\ifcase#1 0\or 1\or 2\or 3\or 4\or 5\or 6\or 7\or 8\or
 9\or A\or B\or C\or D\or E\or F\fi}

\font\teneufm=eufm10
\font\seveneufm=eufm7
\font\fiveeufm=eufm5
\newfam\eufmfam
\textfont\eufmfam=\teneufm
\scriptfont\eufmfam=\seveneufm
\scriptscriptfont\eufmfam=\fiveeufm

\catcode`\@=\csname pre amssym.def at\endcsname

\def    \eps    {\epsilon}

\newcommand{\CM}{{\mathcal M}}

\newcommand{\CS}{{\mathcal S}}

\newcommand{\supp}{\operatorname{supp}}

\newcommand{\id}{{\mathit id}}
\newcommand{\const}{{\mathit const}}

\newcommand{\tF}{\tilde{F}}

\newcommand{\hK}{\hat{K}}
\newcommand{\hH}{\hat{H}}
\newcommand{\tK}{\tilde{K}}
\newcommand{\tf}{\tilde{f}}
\newcommand{\tpsi}{\tilde{\psi}}
\newcommand{\tPhi}{\tilde{\Phi}}
\newcommand{\tH}{\tilde{H}}

\def    \C      {{\mathbb C}}
\def    \R      {{\mathbb R}}

\def    \Z      {{\mathbb Z}}

\def    \12    {{\frac{1}{2}}}

\def    \p      {\partial}

\def    \Sp     {\operatorname{Sp}}
\def    \U     {\operatorname{U}}

\def    \HF     {\operatorname{HF}}
\def    \Ho    {\operatorname{H}}
\def    \HFL     {\operatorname{HF}^{\mathit{loc}}_*}
\def    \HML     {\operatorname{HM}^{\mathit{loc}}_*}
\def    \HFLN     {\operatorname{HF}^{\mathit{loc}}_{n}}
\def    \HMLM     {\operatorname{HM}^{\mathit{loc}}_{m}}
\def    \HM     {\operatorname{HM}}

\def    \CF     {\operatorname{CF}}

\def    \vol     {\operatorname{vol}}

\def    \bx     {\bar{x}}
\def    \by     {\bar{y}}

\def    \MUCZ  {\operatorname{\mu_{\scriptscriptstyle{CZ}}}}

\def    \ssminus        {\smallsetminus}

\begin{document}


\setlength{\smallskipamount}{6pt}
\setlength{\medskipamount}{10pt}
\setlength{\bigskipamount}{16pt}





\title[The Conley Conjecture]{The Conley Conjecture}

\author[Viktor Ginzburg]{Viktor L. Ginzburg}

\address{Department of Mathematics, UC Santa Cruz,
Santa Cruz, CA 95064, USA}
\email{ginzburg@math.ucsc.edu}

\subjclass[2000]{53D40, 37J10, 70H12}
\date{\today}
\thanks{The work is partially supported by the NSF and by the faculty
research funds of the University of California, Santa Cruz.}


\begin{abstract}
We prove the Conley conjecture for a closed symplectically aspherical
symplectic manifold: a Hamiltonian diffeomorphism of a such a manifold
has infinitely many periodic points. More precisely, we show that a
Hamiltonian diffeomorphism with finitely many fixed points has simple
periodic points of arbitrarily large period.  This theorem
generalizes, for instance, a recent result of Hingston establishing
the Conley conjecture for tori.

\end{abstract}

\maketitle

\tableofcontents

\section{Introduction}
\labell{sec:intro}

We show that a Hamiltonian diffeomorphism of a closed symplectically
aspherical manifold has infinitely many periodic points. More
precisely, we prove that such a diffeomorphism with finitely many
fixed points has simple periodic points of arbitrarily large period.
For tori, this fact, recently established by Hingston, \cite{Hi}, was
conjectured by Conley, \cite{Co,SZ} and is frequently referred to as
the \emph{Conley conjecture}.  (See also \cite{FrHa,LeC} and references
therein for similar results for Hamiltonian diffeomorphisms and
homeomorphisms of surfaces.) The proof given here uses some crucial
ideas from~\cite{Hi}, but is completely self-contained.

\subsection{Principal results} The main result of the paper is

\begin{Theorem}
\labell{thm:main} Let $\varphi\colon W\to W$ be a Hamiltonian
diffeomorphism of a closed symplectically aspherical manifold
$W$. Assume that the fixed points of $\varphi$ are isolated. Then
$\varphi$ has simple periodic points of arbitrarily large period.
\end{Theorem}

We refer the reader to Section \ref{sec:general} for the definitions.
Here we only point out that a Hamiltonian diffeomorphism is the
time-one map of a time-dependent Hamiltonian flow and that the
manifolds $W$ with $\pi_2(W)=0$ (e.g., tori and surfaces of genus
greater than zero) are among symplectically aspherical
manifolds. Thus, Theorem \ref{thm:main} implies in particular the
Conley conjecture for tori, \cite{Hi}, and the results of \cite{FrHa}
on Hamiltonian diffeomorphisms of such surfaces.

\begin{Corollary}
\labell{cor:main}
A Hamiltonian diffeomorphism $\varphi$ of a closed symplectically
aspherical manifold has infinitely many simple periodic points.
\end{Corollary}

\begin{Remark}
The example of an irrational rotation of $S^2$ shows that in general
the requirement that $W$ is symplectically aspherical cannot be
completely eliminated; see, however, \cite{FrHa}.
Let $H$ be a periodic in time Hamiltonian giving rise to
$\varphi$. Since periodic points of
$\varphi$ are in one-to-one correspondence with periodic orbits of the
time-dependent Hamiltonian flow $\varphi^t_H$, Theorem \ref{thm:main} and
Corollary \ref{cor:main} can be viewed as results about periodic
orbits of $H$. Then, in both of the statements, the periodic orbits can be
assumed to be contractible. (It is not hard to see that
contractibility is a property of a fixed point rather than of an orbit,
independent of the choice of $H$.) Finally note that, as simple examples
show, the assumption of Theorem \ref{thm:main} that the fixed points of
$\varphi$ are isolated cannot be dropped as long as the periodic orbits
are required to be contractible.
\end{Remark}

There are numerous parallels between the Hamiltonian Conley conjecture
considered here and its Lagrangian counterpart; see, e.g.,
\cite{Lo,Lu,Ma} and references therein. The similarity between the two
problems goes beyond the obvious analogy of the statements and can also easily be
seen on the level of the proofs, although the methods
utilized in \cite{Lo,Lu,Ma} are quite different from the Floer
homological techniques used in the present paper. Thus, for instance,
our Proposition \ref{prop2} plays the same role as Bangert's
homological vanishing method originating from \cite{Ba,BK} in, e.g.,
\cite{Lo,Ma}.

\subsection{Methods}
In the framework of symplectic topology, there are two essentially
different approaches to proving results along the lines of
the Conley conjecture.

The first approach, due to Conley and Salamon and Zehnder,
\cite{CZ86,SZ}, is based on an iteration formula for the
Conley--Zehnder index, asserting that the index of an isolated
weakly non-degenerate orbit either grows linearly under iterations
or its absolute value does not exceed $n-1$, where $2n=\dim W$.
This, in particular, implies that the local Floer homology of such
an orbit eventually becomes zero in degree $n$ as the order of
iteration grows, provided that the orbit remains isolated. (We refer
the reader to Sections \ref{sec:prelim} and \ref{sec:LFH} for the
definitions.  The argument of Salamon and Zehnder, \cite{SZ}, does
not rely on the notion of local Floer homology, but this notion
becomes indispensable in the proof of Theorem \ref{thm:main}.) Since
the Floer homology of $W$ in degree $n$ is non-zero, it follows that
when all one-periodic orbits are weakly non-degenerate, new
simple orbits must be created by large prime iterations to generate
the Floer homology in degree $n$; see \cite{SZ} for details.

The second approach comprises a broad class of methods and is based on
the idea that a Hamiltonian $H$ with sufficiently large variation must
have one-periodic orbits with non-vanishing action.  Since iterating a
Hamiltonian diffeomorphism $\varphi$ has the same effect as, roughly
speaking, increasing the variation of $H$, one can expect $\varphi$ to
have infinitely many periodic points. When a sufficiently accurate
upper or lower bound on the action is available, the orbits can be
shown to be simple. The results obtained along these lines are
numerous and use a variety of symplectic topological techniques and
assumptions on $W$ and $H$.

For instance, if the support of $H$ is displaceable and the variation
of $H$ is greater than the displacement energy $e$ of the support,
one-periodic orbits with action in the range $(0,\,e]$ have been shown
to exist for many classes of symplectic manifolds and Hamiltonians;
see, e.g., \cite{CGK,FH,FHW,FS,Gu,HZ,schwarz,Vi:gen}. Then, the
\emph{a priori} bound on action implies the existence of 
simple periodic orbits with non-zero action and arbitrarily large period. 
These methods do
not rely on particular requirements on the fixed points of $\varphi$,
but the assumption that the support is displaceable appears at this
moment to be crucial. Within this broad class is also a group of
methods applicable to Hamiltonians $H$ with sufficiently degenerate
large or ``flat'' maximum and detecting orbits with action slightly
greater than the maximum of $H$; see, e.g.,
\cite{Gi,GG,Hi,HZ:cap,HZ,Ke:new,KL,LM,MDS,Oh:chain}.  Iterating $\varphi$
can be viewed as stretching $H$ near its maximum, and thus increasing
its variation. Hence, methods from this group can also be used in some
instances to prove the existence of simple periodic orbits of large
period. Here the condition that the maximum is in a certain sense flat
is crucial, but the assumption that the support of the Hamiltonian is
displaceable is less important and, in some cases, not required at
all. In fact, what appears to matter is that the set where the maximum
is attained is relatively small (e.g., symplectic as in \cite{GG} or
displaceable as in \cite{LM} or just isolated as in \cite{Gi,KL}).  It
is one of these methods, combined with the Conley--Salamon--Zehnder
approach, that we use in the proof of the Conley conjecture.

A work of Hingston \cite{Hi} clearly suggests the idea, which is
central to our proof, that the two approaches outlined above can be
extended to cover the case of an arbitrary Hamiltonian. Namely, the
method of \cite{SZ} detects infinitely many simple periodic points of
arbitrarily large periods, unless there exists a strongly degenerate
$\tau$-periodic point $p$ such that the local Floer homology groups of
the $\tau$-th iteration $\varphi^\tau$ and of a large iteration of
$\varphi^\tau$ at $p$ are non-zero in degree $n$; see Section
\ref{sec:proof}. Then we show (Proposition \ref{prop1}) that the
$\varphi^{\tau}$ is the time-$\tau$ flow of a $\tau$-periodic
Hamiltonian, say $H_t$, such that $p$ is a (constant) local maximum of
$H_t$ for all $t$ and this maximum is in a certain sense very
degenerate; cf.\ \cite{Hi}. (However, the Hessian $d^2(H_t)_p$ need
not be identically zero.)  Finally, we prove (Proposition \ref{prop2})
that large iterations of $H$ have periodic orbits with actions
arbitrarily close to the action of the iterated Hamiltonian at $p$.
These orbits are necessarily simple due to the lower and upper bounds
on the action.  Proposition \ref{prop2} is established by using a
simple squeezing argument akin to the ones from \cite{BPS,GG}.
This concludes the proof of the theorem.

This argument is extremely flexible and readily extends to manifolds
convex at infinity or geometrically bounded and wide; see \cite{FS} and
\cite{Gu} for the definitions. We will give a detailed proof of the Conley
conjecture for such manifolds elsewhere.

\subsection{Organization of the paper} In Section \ref{sec:prelim},
we set notation and conventions, briefly review elements of Floer
theory, and also discuss the properties of loops of Hamiltonian
diffeomorphisms relevant to the proof. Local Floer homology is the
subject of Section \ref{sec:LFH}. In Section \ref{sec:proof}, we state
Propositions \ref{prop1} and \ref{prop2} mentioned above and derive
Theorem \ref{thm:main} from these propositions. Proposition
\ref{prop1} reduces the problem to the case of a Hamiltonian with
strict, but ``flat'', local maximum.  This proposition is proved in
Sections \ref{sec:pr-prop1} and \ref{sec:gf} by adapting an argument
from \cite{Hi}. Proposition \ref{prop2} asserting the existence of
simple periodic orbits of large period for such a Hamiltonian and
completing the proof of Theorem \ref{thm:main} is established in
Section~\ref{sec:prop2-pf}. 

\subsection{Acknowledgments} The author is deeply grateful to
Ba\c sak G\"urel, Doris Hein, Ely Kerman, and Felix Schlenk for their numerous
valuable remarks and suggestions.

\section{Preliminaries}
\labell{sec:prelim}

\subsection{Notation and conventions}
\labell{sec:notation}

\subsubsection{General}
\labell{sec:general}

Throughout the paper, a smooth $m$-dimensional manifold is denoted
by $M$ and $W$ stands for a symplectic manifold of dimension $2n$,
equipped with a symplectic form $\omega$. The manifold $(W,\omega)$ is
always assumed to be closed and \emph{symplectically aspherical}, i.e.,
$\omega|_{\pi_2(W)}=0=c_1|_{\pi_2(W)}$, where $c_1$ is the first Chern
class of $W$; see, e.g., \cite{mdsa-book}. An almost complex structure
$J$ on $W$ is said to be \emph{compatible with $\omega$} if
$\omega(\cdot,\, J\cdot)$ is a Riemannian metric on $W$. When $J=J_t$
depends on an extra parameter $t$ (time), this condition is required
to hold for every $t$. A (time-dependent) metric of the form
$\omega(\cdot,\, J\cdot)$ is said to be compatible with $\omega$.

The group of linear symplectic transformations of a finite-dimensional
linear symplectic space $(V,\omega)$ is denoted by $\Sp(V)$. We will also
need the fact that $\pi_1(\Sp(V))\cong\Z$ (see, e.g.,
\cite{mdsa-intro}), and hence $\Ho_1(\Sp(V);\Z)\cong \Z$. To be more specific,
fixing a linear complex structure $J$ on $V$, compatible with
$\omega$, gives rise to an inclusion $\U(V)\hookrightarrow \Sp(V)$ of the
unitary group into the symplectic group. This inclusion is a homotopy
equivalence.  The isomorphism $\pi_1(\Sp(V))\cong\Z$ is the
composition of the isomorphism $\pi_1(\Sp(V))=\pi_1(\U(V))$, the
isomorphism of the fundamental groups induced by $\det\colon U(V)\to
S^1$ (the unit circle in $\C$), and the identification
$\pi_1(S^1)\cong\Z$ arising from fixing the counter clock-wise
orientation of $S^1$. Note that the resulting isomorphism is
independent of the choice of $J$. The \emph{Maslov index} of a loop in
$\Sp(V)$ is the class of this loop in $H_1(\Sp(V);\Z)\cong \Z$.

\subsubsection{Hamiltonians and periodic orbits}
We use the notation $S^1$ for the circle $\R/\Z$ and the circle
$\R/ T\Z$ of circumference $T>0$ is denoted by $S^1_T$.
All Hamiltonians $H$ on $W$ considered in this paper are assumed to be
$T$-periodic (in time), i.e., $H\colon S^1_T\times W\to\R$.
We set $H_t = H(t,\cdot)$ for $t\in S^1_T$. The Hamiltonian vector field
$X_H$ of $H$ is defined by $i_{X_H}\omega=-dH$.

Let $\gamma\colon S^1_T\to W$ be a contractible loop. The action of
$H$ on $\gamma$ is defined by
$$
A_H(\gamma)=A(\gamma)+\int_{S^1_T} H_t(\gamma(t))\,dt.
$$
Here $A(\gamma)$ is the negative symplectic area bounded by $\gamma$, i.e.,
$$
A(\gamma)=-\int_z\omega,
$$
where $z\colon D^2\to W$ is such that $z|_{S^1_T}=\gamma$.

The least action principle asserts that the critical points of $A_H$ on
the space of all contractible maps $\gamma\colon S^1_T\to W$
are exactly the contractible $T$-periodic orbits of the time-dependent
Hamiltonian flow $\varphi_H^t$ of $H$.
When the period $T$ is clear from the context and, in particular, if
$T=1$, we denote the time-$T$ map
$\varphi^T_H$ by $\varphi_H$. The action spectrum
$\CS(H)$ of $H$ is the set of critical values of $A_H$.
This is a zero measure, closed set; see, e.g., \cite{HZ,schwarz}.
In this paper we are only concerned with contractible periodic orbits.
\emph{A periodic orbit is always assumed to be contractible, even if
this is not explicitly stated}.

\begin{Definition}
A $T$-periodic orbit $\gamma$ of $H$ is \emph{non-degenerate} if
the linearized return map
$d\varphi^T_H \colon T_{\gamma(0)}W\to T_{\gamma(0)}W$ has no eigenvalues
equal to one. Following \cite{SZ}, we call $\gamma$
\emph{weakly non-degenerate} if at least one of the eigenvalues is different
from one. When all eigenvalues are equal to one, the orbit is said to be
\emph{strongly degenerate}.
\end{Definition}

When $\gamma$ is non-degenerate or even weakly non-degenerate, the
so-called \emph{Conley--Zehnder index} $\MUCZ(\gamma)\in\Z$ is
defined, up to a sign, as in \cite{Sa,SZ}. More specifically, in this
paper, the Conley--Zehnder index is the negative of that of
\cite{Sa}. In other words, we normalize $\MUCZ$ so that
$\MUCZ(\gamma)=n$ when $\gamma$ is a non-degenerate maximum of an
autonomous Hamiltonian with small Hessian. More generally, when $H$ is
autonomous and $\gamma$ is a non-degenerate critical point of $H$ such
that the eigenvalues of the Hessian (with respect to a metric
compatible with $\omega$) are less than $2\pi/ T$, the Conley--Zehnder
index of $\gamma$ is equal to one half of the negative signature of
the Hessian.

Sometimes, the same Hamiltonian can be treated as
$T$-periodic for different values of $T>0$. For instance, an
autonomous Hamiltonian is $T$-periodic for every $T>0$ and a
$T$-periodic Hamiltonian can also be viewed as $Tk$-periodic for
any integer $k$. In this paper, it will be essential to keep track
of the period. Unless specified otherwise, every Hamiltonian $H$
considered here is originally one-periodic and $T$ is always an
integer. When we wish to view $H$ as a $T$-periodic Hamiltonian, we
denote it by $H^{(T)}$ and refer to it as the \emph{$T$-th iteration}
of $H$. (The parentheses here are used to distinguish iterated
Hamiltonians from families of Hamiltonians, say $H^s$, parametrized
by $s$.)  Since $H^{(T)}$ is regarded as a $T$-periodic Hamiltonian,
it makes sense to speak only about $T$-periodic (or $Tk$-periodic)
orbits of $H^{(T)}$. Clearly, $T$-periodic orbits of $H^{(T)}$ are simply
$T$-periodic orbits of $H$.  When $\gamma\colon S^1\to W$ is a
one-periodic orbit of $H$, its \emph{$T$-th iteration} is the
obviously defined map $\gamma^{(T)}\colon S^1_T\to W$ obtained by
composing the $T$-fold covering map $S^1_T\to S^1$ with $\gamma$.
Thus, $\gamma^{(T)}$ is a $T$-periodic orbit of $H$ and
$H^{(T)}$.  We call a $T$-periodic orbit \emph{simple} if it is not
an iteration of an orbit of a smaller period.

As is well-known, the fixed points of $\varphi_H :=\varphi_H^1$ are in
one-to-one correspondence with (not-necessarily contractible)
one-periodic orbits of $H$. Likewise, the $T$-periodic points of
$\varphi_H$, i.e., the fixed points of $\varphi_H^T$, are in
one-to-one correspondence with (not-necessarily contractible)
$T$-periodic orbits. In the proof of Theorem \ref{thm:main}, we will
work with (contractible!) periodic orbits of a Hamiltonian $H$ whose
time-one flow is $\varphi$. In fact, as is easy to see from Section
\ref{sec:loops}, the free homotopy type of the one-periodic orbit of
$H$ through a fixed point $p$ of $\varphi$ is completely determined by
$\varphi$ and $p$ and is independent of the choice of $H$. The same holds
for $T$-periodic points and orbits. Hence, ``contractible fixed points
or periodic points'' of $\varphi$, i.e., those with contractible
orbits, are well defined and we will establish Theorem \ref{thm:main}
for points in this class.

When $K$ and $H$ are two (say, one-periodic) Hamiltonians, the composition
$K\# H$ is defined by the formula
$$
(K\# H)_t=K_t+H_t\circ(\varphi^t_K)^{-1}.
$$
The flow of $K\# H$ is $\varphi^t_K\circ \varphi^t_H$. In general, $K\# H$
is not one-periodic. However, this will be the case if, for example,
$H_0\equiv 0\equiv H_1$.
Another instance when the composition $K\# H$ of two one-periodic
Hamiltonians is automatically one-periodic is when the flow
$\varphi^t_K$ is a loop of Hamiltonian diffeomorphisms, i.e.,
$\varphi_K^1=\id$.

\subsubsection{Norms with respect to a coordinate system}
\labell{sec:norms}

In what follows, it will be convenient to use
a somewhat unconventional terminology and work with $C^k$-norms of
functions, vector fields, etc. taken with respect to a coordinate
system.

Let $\xi$ be a coordinate system on a neighborhood $U$ of a point
$p\in M^m$, i.e., $\xi$ is a diffeomorphism $U\to \xi(U)\subset \R^m$
sending $p$ to the origin. (Thus, $U$ is a part of the data $\xi$.) Let $f$
be a function on $U$ or on the entire manifold $M$. The $C^k$-norm
$\|f\|_{C^k(\xi)}$ of $f$ with respect to $\xi$
is, by definition, the $C^k$-norm of $f$ on $U$ with
respect to the flat metric associated with $\xi$, i.e., the pull-back by $\xi$
of the standard metric on $\R^m$. The $C^k$-norm with respect to $\xi$
of a vector field or a field of operators on $U$ is defined in a
similar fashion.

Likewise, the norm $\|v\|_\Xi$ of a vector $v$ in a finite-dimensional
vector space $V$ with respect to a basis $\Xi$ is the norm of
$v$ with respect to the Euclidean inner product for which $\Xi$ is an
orthonormal basis. The norm of an operator $V\to V$ with respect to $\Xi$
is defined in a similar way. When $\xi$ is a coordinate
system near $p$, we denote by $\xi_p$ the natural coordinate
basis in $T_pM$ arising from $\xi$.

\begin{Example}
\labell{exam:norm1} Let $A\colon V\to V$ be a linear map with all
eigenvalues equal to zero.  Then $\|A\|_{\Xi}$ can be made
arbitrarily small by choosing a suitable basis $\Xi$. In other
words, for any $\sigma>0$ there exists $\Xi$ such that
$\|A\|_{\Xi}<\sigma$. Indeed, in some basis $A$ is given by an upper
triangular matrix with zeros on the diagonal; $\Xi$ is then obtained
by appropriately scaling the elements of this basis.
\end{Example}

\begin{Example}
\labell{exam:norm2} Restricting $\xi$ to a smaller neighborhood of
$p$ reduces the norm of $f$. However, one cannot make,
say, $\|f\|_{C^1(\xi)}$ arbitrarily small by shrinking $U$ unless
$f(p)=0$ and $df_p=0$. Indeed, $\|f\|_{C^1(\xi)}\geq\max\{|f(p)|,
\|df_p\|_{\xi_p}\}$. It is clear that for a fixed basis $\Xi$ in
$T_pM$ and a function $f$ near $p$ there exists a coordinate system
$\xi$ with $\xi_p=\Xi$, such that $\|f\|_{C^1(\xi)}$ is arbitrarily
close to $\max\{|f(p)|, \|df_p\|_{\Xi}\}$.
\end{Example}

\subsection{Floer homology}

In this section, we briefly recall the notion of filtered Floer
homology for closed symplectically aspherical manifolds.  All
definitions and results mentioned here are quite standard and well
known and we refer the reader to Floer's papers
\cite{F:Morse,F:grad,F:c-l,F:witten}, to
\cite{BPS,FHW,FH,FHS,SZ,schwarz}, or to
\cite{HZ,mdsa-book,Sa:london,Sa} for introductory accounts of the
construction of Floer homology in this (or more general) setting.

\subsubsection{Definitions}
Let us first focus on one-periodic Hamiltonians. Consider a
Hamiltonian $H$ such that all one-periodic orbits of $H$ are
non-degenerate. This is a
generic condition and we will call such Hamiltonians
\emph{non-degenerate}. Let $J=J_t$ be a (time-dependent: $t\in S^1$)
almost complex structure on $W$ compatible with $\omega$. For two
one-periodic orbits $x^\pm$ of $H$ denote by $\CM_H(x^-,x^+,J)$ the
space of solutions $u\colon S^1\times \R\to W$ of the Floer equation
\begin{equation}
\label{eq:floer}
\frac{\p u}{\p s}+ J_t(u) \frac{\p u}{\p t}=-\nabla H_t(u)
\end{equation}
which are asymptotic to $x^\pm$ at $\pm\infty$, i.e.,
$u(s,t)\to x^\pm(t)$ point-wise as $s\to \pm\infty$.

The energy $E(u)$ of a solution $u$ of the Floer equation, \eqref{eq:floer},
is
$$
E(u)=\int_{-\infty}^\infty \left\|\frac{\p u}{\p s}\right\|_{L^2(S^1)}^2\,ds
=\int_{-\infty}^\infty \int_{S^1}\left\|\frac{\p u}{\p t}-J\nabla H (u)
\right\|^2 \,dt\,ds,
$$
where we set $u(s):=u(s,\,\cdot)\colon S^1\to W$. Every finite energy solution
of \eqref{eq:floer} is asymptotic to some $x^\pm$ and
$$
A_H(x^-)-A_H(x^+)=E(u).
$$

When $J$ meets certain standard regularity requirements that hold
generically (see, e.g., \cite{FHS,SZ}), the space $\CM_H(x^-,x^+,J)$
is a smooth manifold of dimension $\MUCZ(x^+)-\MUCZ(x^-)$.  This
space carries a natural $\R$-action $(\tau\cdot u)(t,s)=u(t,s+\tau)$
and we denote by $\hat{\CM}_H(x^-,x^+,J)$ the quotient
$\CM_H(x^-,x^+,J)/\R$. When $\MUCZ(x^+)-\MUCZ(x^-)=1$, the set
$\hat{\CM}_H(x^-,x^+,J)$ is finite and we denote the number,
$\!\!\!\mod 2$, of
points in this set by $\#\big(\hat{\CM}_H(x,y,J)\big)$.

Let $a<b$ be two points outside $\CS(H)$. Denote by
$\CF_k^{(a,\,b)}(H)$ the vector space over $\Z_2$ generated by
one-periodic orbits of $H$ with $\MUCZ(x)=k$ and $a<A_H(x)<b$. The Floer
differential
$$
\p \colon \CF_k^{(a,\,b)}(H)\to \CF_{k-1}^{(a,\,b)}(H)
$$
is defined by
$$
\p x=\sum_y \#\big(\hat{\CM}_H(x,y,J)\big)\cdot y,
$$
where the summation is over all $y$ such that $\MUCZ(y)=k-1$ and
$a<A_H(y)<b$.  As is well-known, $\p^2=0$.  The homology
$\HF_*^{(a,\,b)}(H)$ of the resulting complex is called the
\emph{filtered Floer homology} of $H$ for the interval $(a,\,b)$.
Thus, $\HF_*(H):=\HF^{(-\infty,\infty)}_*(H)$ is the ordinary
Floer homology. It is a standard fact that $\HF_*(H)=\Ho_{*+n}(W;\Z_2)$.  In
general, $\HF^{(a,\,b)}_*(H)$ depends on the Hamiltonian $H$, but not
on~$J$.

The subcomplexes $\CF_*^{(-\infty,\,b)}(H)$, where $b\in \R\ssminus
\CS(H)$, form a filtration of the total Floer complex
$\CF_*(H):=\CF_*^{(-\infty,\,\infty)}(H)$, called the \emph{action
filtration}, and $\CF_*^{(a,\,b)}(H)$ can be identified with the
complex $\CF_*^{(-\infty,\,b)}(H)/\CF_*^{(-\infty,\,a)}(H)$;
see, e.g., \cite{FH,schwarz}.  Let now $a<b<c$ be three points outside
$\CS(H)$.  Then, similarly, $\CF^{(a,\,b)}_*(H)$ is a subcomplex of
$\CF^{(a,\,c)}_*(H)$, and $\CF^{(b,\,c)}_*(H)$ is naturally isomorphic
to $\CF^{(a,\,c)}_*(H)/\CF^{(a,\,b)}_*(H)$. As a result, we have the
long exact sequence
\begin{equation}
\label{eq:seq}
\ldots\to \HF^{(a,\,b)}_*(H)\to\HF^{(a,\,c)}_*(H)\to \HF^{(b,\,c)}_*(H)\to
\HF^{(a,\,b)}_{*-1}(H)\to \ldots .
\end{equation}

The filtered Floer homology of $H$ is defined even when the periodic
orbits of $H$ are not necessarily non-degenerate, provided that $a<b$ are
outside $\CS(H)$.  Namely, let $\tH$ be a $C^2$-small
perturbation\footnote{For the sake of brevity, we refer to the
function $\tH$, rather than the difference $\tH-H$, as a perturbation
of $H$. However, it is the difference $\tH-H$ that is required to be
$C^2$-small.} of $H$ with non-degenerate one-periodic orbits.  The
filtered Floer homology $\HF^{(a,\,b)}_*(H)$ of $H$ is by definition
$\HF^{(a,\,b)}_*(\tH)$. (Clearly, $a<b$ are still outside $\CS(\tH)$.)
These groups are canonically isomorphic for different choices of $\tH$
(close to $H$) and the results discussed here hold for
$\HF^{(a,\,b)}_*(H)$, \cite{BPS,FH,FHW,schwarz,vi:functors}. In fact,
when an assertion concerns individual Hamiltonians (as opposed to
families of Hamiltonians), it is usually sufficient to prove the
assertion in the non-degenerate case, for then it extends ``by
continuity'' to all Hamiltonians.

These constructions carry over to $T$-periodic Hamiltonians
word-for-word by replacing one-periodic orbits with $T$-periodic ones. When
$H$ is one-periodic, but we treat it as $T$-periodic for some integer
$T>0$, we denote the resulting Floer homology groups
$\HF^{(a,\,b)}_*\big(H^{(T)}\big)$.

\subsubsection{Homotopy maps}
\labell{sec:homotopy} Consider two non-degenerate Hamiltonians $H^0$
and $H^1$.  Let $H^s$ be a homotopy from $H^0$ to $H^1$. By
definition, this is a family of Hamiltonians parametrized by $s\in \R$
such that $H^s\equiv H^0$ when $s$ is large negative and $H^s\equiv
H^1$ when $s$ is large positive.  (Strictly speaking, the notion of
homotopy includes also a family of almost complex structures $J_s$;
see, e.g., \cite{BPS,FH,FHW,Sa}. We suppress this part of the homotopy
structure in the notation.)  Assume, in addition, that the homotopy is
\emph{monotone decreasing}, i.e., $H^s_t(p)$ is a decreasing function of $s$
for all $p\in W$ and $t\in S^1$. (Thus, in particular, $H^0\geq H^1$.)
Then, whenever $a<b$ are outside $\CS(H^0)$ and $\CS(H^1)$, the
homotopy $H^s$ induces a homomorphism of complexes
$\Psi_{H^0,H^1}\colon\CF_*^{(a,\,b)}(H^0)\to \CF_*^{(a,\,b)}(H^1)$ by
the standard continuation construction; see, e.g.,
\cite{BPS,FHW,FH,schwarz}. Namely, for a one-periodic orbit $x$ of
$H^0$ and a one-periodic orbit $y$ of $H^1$, let $\CM_H(x,y,J)$ be the
space of solutions of \eqref{eq:floer} with $H^s$ on the right hand
side, asymptotic to $x$ and, respectively, $y$ at $\pm\infty$.  Under
the well-known regularity requirements on $J$ and $H^s$, the space
$\CM_{H^s}(x,y,J)$ is a smooth manifold of dimension
$\MUCZ(y)-\MUCZ(x)$; see, e.g.,
\cite{FH,Sa,SZ,schwarz,schwarz:book}. Moreover, $\CM_{H^s}(x,y,J)$ is
a finite collection of points when $\MUCZ(y)=\MUCZ(x)$.  The map
$\Psi_{H^0,H^1}$ is defined by
\begin{equation}
\labell{eq:homo}
\Psi_{H^0,H^1}(x)
=\sum_y \#\big({\CM}_H(x,y,J)\big)\cdot y,
\end{equation}
where the summation is over all $y$ such that $\MUCZ(y)=\MUCZ(x)$ and
$a<A_H(y)<b$.

The induced \emph{homotopy map} in the filtered Floer homology, also
denoted by $\Psi_{H^0,H^1}$, is independent of the decreasing homotopy
$H^s$ and commutes
with the maps from the long exact sequence \eqref{eq:seq}, see,
e.g., \cite{BPS,FH,Sa,SZ,schwarz:book,schwarz}. By ``continuity'' in
the Hamiltonians and the homotopy, this construction extends to all
(not necessarily non-degenerate) Hamiltonians and all decreasing
(but not necessarily regular) homotopies as long as $a<b$ are not in
$\CS(H^0)$ and $\CS(H^1)$.

A (non-monotone) homotopy $H^s$ from $H^0$ to $H^1$ with $a$ and $b$
outside $\CS(H^s)$ for all $s$ gives rise to an isomorphism between
the groups $\HF^{(a,\,b)}_*(H^s)$, and hence, in particular,
\begin{equation}
\labell{eq:viterbo}
\HF^{(a,\,b)}_*(H^0)\cong \HF^{(a,\,b)}_*(H^1);
\end{equation}
see \cite{BPS,vi:functors}. This isomorphism is defined by breaking
the homotopy $H^s$ into a composition of homotopies $H^{i,s}$ close to
the identity. Each of the homotopies $H^{i,s}$ and its inverse homotopy
increase action by no more than some small $\eps>0$. Then, it is shown
that the map in $\HF^{(a,\,b)}_*$ induced by $H^{i,s}$
is an isomorphism. Although this construction requires
additional choices, it is not hard to see that the isomorphism
\eqref{eq:viterbo} is uniquely determined by the homotopy $H^s$.
Furthermore, \eqref{eq:viterbo} commutes with the maps from the long
exact sequence \eqref{eq:seq}, provided that all three points $a<b<c$
are outside $\CS(H^s)$ for all $s$; \cite{Gi}. Note also that when
$H^s$ is a decreasing homotopy, the isomorphism  \eqref{eq:viterbo}
coincides with $\Psi_{H^0,H^1}$.

\begin{Example}
\labell{ex:isospec} A homotopy $H^s$ is said to be
\emph{isospectral} if $\CS(H^s)$ is independent of $s$. In this case,
the isomorphism \eqref{eq:viterbo} is defined for any $a<b$ outside
$\CS(H^s)$.

For instance, let $\psi_s^t$, where $t\in S^1$ and $s\in [0,\, 1]$, be
a family of loops of Hamiltonian diffeomorphisms based at $\id$, i.e.,
$\psi^0_s=\id$ for all $s$. In other words, $\psi_s^t$ is a based
homotopy from the loop $\psi_0^t$ to the loop $\psi_1^t$. Let $G^s_t$
be a family of one-periodic Hamiltonians generating these loops and
let $H$ be a fixed one-periodic Hamiltonian. Then $H^s:=G^s\# H$ is an
isospectral homotopy, provided that $G^s$ are suitably normalized.
(Namely, $A(G^s)=0$ for all $s$; see Section \ref{sec:loops}.)
\end{Example}

It is easy to see that if $K\geq H^s$ for all $s$,
the isomorphism \eqref{eq:viterbo} intertwines monotone homotopy
homomorphisms from $K$ to $H^0$ and to $H^1$, i.e.,
the diagram
\begin{equation}
\labell{eq:diag-iso}
\xymatrix{
{\HF^{(a,\,b)}_*(K)} \ar[d]_{\Psi_{K,H^0}}\ar[rd]^{\Psi_{K,H^1}} &\\
{\HF^{(a,\,b)}_*(H^0)} \ar[r]^{\cong} & {\HF^{(a,\,b)}_*(H^1)}
}
\end{equation}
is commutative. Note that it is not at all clear whether the same is
true if we only require that $K\geq H^0$ and $K\geq H^1$.

\subsubsection{Non-triviality criterion for a homotopy map}
We conclude this section by establishing a technical result, used
later on in the proof, giving a criterion for a monotone homotopy map
to be non-zero.

\begin{Lemma}
\labell{lemma:non-zero} Let $H^s$ be a monotone decreasing homotopy
such that a point $p$ is a non-degenerate constant one-periodic orbit
of $H^s$ and $H^s_t(p)=c$ for all $s$ and $t$.  Then the monotone
homotopy map
$$
\Psi_{H^0,H^1}\colon \HF^{(a,\,b)}_*(H^0)
\to\HF^{(a,\,b)}_*(H^1)
$$
is non-trivial, provided that $\CS(H^0)\cap (a,\,b)=\{c\}
=\CS(H^1)\cap (a,\,b)$.
\end{Lemma}

\begin{Remark}
\labell{rmk:non-zero} In fact, we will prove a stronger result. Let us
perturb $H^0$ and $H^1$ away from $p$, making these Hamiltonians
non-degenerate. Then $p$ is a cycle in $\CF^{(a,\,b)}_*(H^0)$
and $\CF^{(a,\,b)}_*(H^1)$ and, moreover, this
cycle is not homologous to any cycle that does not
include $p$.  (This is easy to see from the energy estimates;
see, e.g., \cite{Sa:london,Sa}.)  In particular, $[p]\neq 0$ in
$\HF^{(a,\,b)}_*(H^0)$ and $\HF^{(a,\,b)}_*(H^1)$ and we will show
that $\Psi_{H^0,H^1}$ sends $[p]\in\HF^{(a,\,b)}_*(H^0)$ to
$[p]\in \HF^{(a,\,b)}_*(H^1)$.

Moreover, a simple modification of our argument proves the following:
There exist $C^2$-small non-degenerate perturbations $\hH^0$
of $H^0$ and $\hH^1$ of $H^1$ for which $p$ is still a non-degenerate
constant one-periodic orbit, and a regular monotone decreasing homotopy
$\hH^s$ from $\hH^0$ to $\hH^1$ such that the cycle $p$ for $\hH^0$
is connected to  $p$  for $\hH^1$ by an odd number of homotopy trajectories
and all such trajectories are contained in a small neighborhood of $p$.
(Note that we do not assume that $p$ is a  non-degenerate
constant one-periodic orbit of $\hH^s$ for all $s$ or that
$\hH^0(p)=\hH^1(p)$.
Of course, the lemma can be further generalized. For instance, the
constant orbit $p$ can be replaced by a fixed one-periodic orbit.
\end{Remark}

\begin{proof}[Proof of Lemma \ref{lemma:non-zero}]
Let us perturb the homotopy on the complement of a neighborhood $U$ of
$p$, keeping the homotopy monotone decreasing, to ensure that all but
a finite number of the Hamiltonians $H^s$ are non-degenerate. In
particular, we will assume that $H^0$ and $H^1$ are such. This can be
achieved by an arbitrarily small perturbation of $H^s$.  We keep the notation
$H^s$ for the perturbed homotopy.

If the homotopy $H^s$ were regular, we would simply argue that the
constant connecting trajectory $u\equiv p$ is the only connecting
trajectory from $p$ for $H^0$ to $p$ for $H^1$. Indeed, $0\leq
E(v)=A_{H^0}(p)-A_{H^1}(p)=0$ for any such connecting trajectory $v$,
and thus $v$ must be constant.  However, while it is easy to guarantee
that $H^s$ is regular away from $p$ by reasoning as in, e.g.,
\cite{FH,FHS}, it is not \emph{a priori} obvious that the
transversality requirements can be satisfied for $u\equiv p$ because of the
constraint $H^s_t(p)=c$. Rather than checking regularity of $u$ by a
direct calculation, we chose to circumvent this difficulty.

As in the proof of \eqref{eq:viterbo} (see, e.g.,
\cite{BPS,Gi,vi:functors}), we can break
the homotopy $H^s$ by reparametrization of $s$
into a composition of homotopies $K^{i,s}$ from
$K^i:=H^{s_i}$ to $K^{i+1}:=H^{s_{i+1}}$ with $K^0=H^0$ and
$K^k=H^1$. These homotopies are monotone, since $H^s$ is monotone,
and close to the identity homotopy. For every $\eps>0$, this can
be done so that the inverse homotopy $\Psi_{K^{i+1},K^{i}}$ from
$K^{i+1}$ to $K^i$ increases the actions by no more than $\eps>0$.
Without loss of generality, we may assume that the
Hamiltonians $K^i$ are non-degenerate. Since all ``direct''
homotopies are monotone decreasing, we have
$$
\Psi_{H^0,H^1}=\Psi_{K^{k-1},K^{k}}\circ\ldots\circ\Psi_{K^1,K^2}
\circ\Psi_{K^0,K^1}.
$$

Observe that it suffices to establish the lemma for $a$ and $b$
arbitrarily close to $c$.  Let $U$ be so small that $p$ is the only
one-periodic trajectory of $H^s$ entering $U$ for all $s$.  (Since $p$
is isolated for all $H^s$, such a neighborhood $U$ exists.) There
exists a constant $\eps_U>0$ such that every Floer anti-gradient
trajectory $v$ connecting $p$ with any other one-periodic orbit with
action in the range $(a,\,b)$ has energy $E(v)>\eps_U$ for any regular
Hamiltonian in the family $H^s$ (cf.\ \cite{Sa:london,Sa}). In
particular, this holds for $K^i$ and $K^{i+1}$ and, moreover, for
every regular Hamiltonian in the homotopy $K^{i,s}$. We will pick $a$
and $b$ so that $c-\eps_U<a$ and $b< c+\eps_U$. Then, for every $K^i$
the point $p$ is a cycle in $\CF_*^{(a,\,b)}(K^i)$ and $p$ is not
homologous to any cycle that does not include $p$.

Now it is sufficient to prove that $\Psi_{K^i,K^{i+1}}$ sends
$[p]\in\HF^{(a,\,b)}_*(K^i)$ to $[p]\in\HF^{(a,\,b)}_*(K^{i+1})$. To
this end, let us first prove that $\Psi_{K^i,K^{i+1}}([p])\neq 0$.  We
may assume that none of the points $a$, $a+\eps$, $b$, $b+\eps$ is
in $\CS(K^i)$ or in $\CS(K^{i+1})$.

It is easy to see (see, e.g.,
\cite{Gi}) that
\begin{equation}
\labell{eq:qi}
\Psi_{K^{i+1},K^{i}}\circ\Psi_{K^i,K^{i+1}}
\colon
\HF^{(a,\,b)}_*(K^i)
\to\HF^{(a+\eps,\,b+\eps)}_*(K^i)
\end{equation}
is the natural ``quotient-inclusion'' map, i.e., the composition of
the ``quotient'' and ``inclusion'' maps $ \HF^{(a,\,b)}_*(K^i)\to
\HF^{(a+\eps,\,b)}_*(K^i) \to\HF^{(a+\eps,\,b+\eps)}_*(K^i).  $
Note that
$\eps_U$ is completely determined by $H^s$ and $U$ and is
independent of how $H^s$
is broken into the homotopies $K^{i,s}$, and thus of $\eps>0$.
Pick $\eps>0$ so that $a+\eps<c$ and $b+\eps<c+\eps_U$. Then $p$ is a
cycle in $\CF_*^{(a+\eps,\,b+\eps)}(K^i)$, which is not homologous to
any cycle that does not include $p$. As a consequence, $[p]\neq 0$
in both of the Floer homology groups in \eqref{eq:qi} and
$\Psi_{K^{i+1},K^{i}}\circ\Psi_{K^i,K^{i+1}}([p])=[p]$.
Therefore, $\Psi_{K^i,K^{i+1}}([p])\neq 0$ in
$\HF^{(a,\,b)}_*(K^{i+1})$.

To show that $\Psi_{K^i,K^{i+1}}([p])=[p]$, we need to refine our
choice of $\eps_U$. Note that there exists $\eps_U>0$ such that, in
addition to the above requirements, every $K^{i,s}$-homotopy
trajectory starting at $p$ and leaving $U$ has energy greater than
$\eps_U$. Then, clearly, the same is true for every sufficiently
$C^2$-small perturbation $\hK^{i,s}$ of $K^{i,s}$.  Again, $\eps_U$
depends only on $H^s$ and $U$, but not on breaking $H^s$ into the
homotopies $K^{i,s}$. (The existence of $\eps_U>0$ with these
properties readily follows from energy estimates for connecting
trajectories, cf.\ \cite{Sa:london,Sa}.)  Pick a $C^2$-small regular
perturbation $\hK^{i,s}$ of $K^{i,s}$. We may still assume that
$\hK^{i,s}$ is monotone decreasing and $p$ is a
non-degenerate constant orbit of $\hK^{i}$ and $\hK^{i+1}$. However,
$p$ is not required to be a constant one-periodic orbit of
$\hK^{i,s}$ for all $s$; nor is $\hK^{i,s}_t(p)$ constant as a function of
$s$ and $t$.  Clearly, the inverse homotopy
to $\hK^{i,s}$ does not increase action by more than $\eps>0$. Hence,
in the homological analysis of the previous paragraph we can replace
$K^i$ and $K^{i+1}$ by $\hK^i$ and $\hK^{i+1}$. In fact, the original
and perturbed Hamiltonians have equal filtered Floer homology for
relevant action intervals and the maps $\Psi_{K^i,K^{i+1}}$ and
$\Psi_{K^{i+1},K^{i}}$ are induced by the maps for $\hK^{i,s}$ acting
on the level of complexes.  Therefore,
$\Psi_{\hK^i,\hK^{i+1}}([p])=\Psi_{K^i,K^{i+1}}([p])\neq 0$ in
$\HF^{(a,\,b)}_*(\hK^{i+1})=\HF^{(a,\,b)}_*(K^{i+1})$, and thus
$\Psi_{\hK^i,\hK^{i+1}}(p)\neq 0$ in $\CF^{(a,\,b)}_*(\hK^{i+1})$.
Since $c-\eps_U<a$ and $c+\eps_U<b$ and every connecting orbit leaving
$U$ must have energy greater than $\eps_U$, we conclude that
$\Psi_{\hK^i,\hK^{i+1}}(p)=p$ in the Floer complexes, and hence in
the Floer homology.
\end{proof}

\subsection{Loops of Hamiltonian diffeomorphisms}
\labell{sec:loops}

In this section, we recall a few well-known facts about loops of
Hamiltonian diffeomorphisms of $W$. We will focus on loops
parametrized by $S^1$, but obviously all results discussed here hold
for loops of any period. Furthermore, throughout the paper all loops
$\psi=\psi^t$ are assumed to be based at $\id$, i.e., $\psi^1=\psi^0=\id$;
contractible loops are thus required to be contractible in this class.

Recall that, as is proved in \cite{schwarz}, the filtered Floer
homology of the Hamiltonian $H$ is determined, up to a shift of
filtration, entirely by the time-one map $\varphi_H$ and is
independent of the Hamiltonian $H$.  This fact translates to geometric
properties of loops of Hamiltonian diffeomorphisms, which are briefly
reviewed below, and is actually proved by first establishing these
properties.

\subsubsection{Global loops}
Let $\psi^t=\varphi^t_G$, $t\in S^1$, be a loop generated by a
periodic Hamiltonian $G$. Then all
orbits $\gamma(t)=\psi^t(p)$ of $\psi^t$ with $t\in S^1$ and $p\in W$ are
one-periodic and lie in the same homotopy class. Hence, every orbit
of $G$ is contractible by the Arnold conjecture.
The action $A(G):=A_G(\gamma)$ is independent of $p\in W$ (see,
e.g., \cite{HZ, schwarz}) and $A(G)=\vol (W)^{-1}\int_0^1\int_W
G\,\omega^n\, dt$, where $\vol(W)$ is the symplectic
volume of $W$.  The latter identity is easy to prove when the loop
$\psi^t$ is contractible (see, e.g., \cite{Gi, schwarz}); in the general case,
this is a non-trivial result, \cite{schwarz}.

For $\gamma$ as above pick a trivialization of $TW|_\gamma$ that
extends to a trivialization of $TW$ along a disk bounded by
$\gamma$. Using this trivialization, we can view the maps
$d\psi^t\colon T_{\gamma(0)}W\to T_{\gamma(t)}W$ as a loop in
$\Sp(T_{\gamma(0)}W)$. Hence, the linearization
$d\psi^t$ along $\gamma$ gives rise to an element in
$\pi_1(\Sp(T_{\gamma(0)}W))=\Z$, which could be called the Maslov
index, $\mu(\psi)$, of the loop $\psi^t$ if it were non-trivial; cf.\
\cite{SZ,Sa}. The Maslov index is well defined: it is independent of
$\gamma$, the trivialization and the disc. (The latter follows from
the fact that $c_1(W)|_{\pi_2(W)}=0$.)  However, as is well known and
as we will soon reprove, $\mu(\psi)=0$; see also, e.g.,
\cite{schwarz}.

Let $H$ be a periodic Hamiltonian on $W$. Recall that $G\#H$ is the
Hamiltonian generating the flow $\varphi^t_G\varphi^t_H$. This
Hamiltonian is automatically one-periodic and its time-one map is
$\varphi_H$. We claim that there exists an isomorphism of filtered
Floer homology
\begin{equation}
\labell{eq:FH-loops}
\HF^{(a,\,b)}_*(H)\cong \HF^{(a+A(G),\,b+A(G))}_*(G\# H).
\end{equation}
Indeed, composition with $\psi^t=\varphi^t_G$ sends one-periodic
orbits of $H$ to one-periodic orbits of $G\# H$ with shift of action
by $A(G)$ and shift of Conley--Zehnder indices by $-2\mu(\psi)$. (See,
e.g., \cite{Sa, schwarz}; the negative sign is a result of the
difference in conventions.) Furthermore, let $u$ be a Floer
anti-gradient trajectory for $H$ and a time-dependent almost complex
structure $J$. Then, as a straightforward calculation shows,
$\psi(u)(s,t):=\psi^t(u(s,t))$ is a Floer anti-gradient trajectory for
$G\# H$ and the almost complex structure $\tilde{J}_t:=d\psi^t \circ
J_t \circ (d\psi^t)^{-1}$. Furthermore, it is clear that the
transversality requirements are satisfied for $(H, J)$ if and only if
they are satisfied for $(G\# H,\tilde{J})$. Therefore, the composition
with $\psi$ commutes with the Floer differential and thus induces an
isomorphism of Floer complexes (and hence homology groups) shifting
action by $A(G)$ and grading by $-2\mu(\psi)$.  Applying this
construction to the full Floer homology $\HF_*(H)\cong H_{*+n}(W)\cong
\HF_*(G\# H)$, we see that the grading shift must be zero, i.e.,
$\mu(\psi)=0$.

\begin{Remark}
When the loop $\psi$ is contractible, the existence of an isomorphism
\eqref{eq:FH-loops} readily follows from \eqref{eq:viterbo}; see
Example \ref{ex:isospec}. However, it is not clear whether this is the
same isomorphism as constructed above.
\end{Remark}

\subsubsection{Local loops}
Let now $\psi^t$ be a loop of (the germs of) Hamiltonian diffeomorphisms
at $p\in W$ generated by $G$. In other words, the maps $\psi^t$ and
the Hamiltonian $G$ are defined on a small neighborhood of $p$ and
$\psi^t(p)=p$ for all $t\in S^1$. Then the action $A(G)$ and the Maslov
index $\mu(\psi)$ are introduced exactly as above with the orbit $\gamma$
taken sufficiently close to $p$. In fact, we can set
$\gamma\equiv p$ and hence $A(G)=\int_0^1 G_t(p)\, dt$ and $\mu(\psi)$ is
just the Maslov index of the loop $d\psi^t_p$ in $\Sp(T_pW)$. Note that in
this case $\mu(\psi)$ need not be zero.

We conclude this section by giving a necessary and sufficient
condition, to be used later, for $\psi$ to extend to a loop of global
Hamiltonian diffeomorphisms of $W$.

\begin{Lemma}
\labell{lemma:loop-ext}
Let $\psi^t$, $t\in S^1$, be a loop of germs of Hamiltonian diffeomorphisms
at $p\in W$. The following conditions are equivalent:

\begin{itemize}
\item [(i)] the loop $\psi$ extends to a loop of global
Hamiltonian diffeomorphisms of $W$,

\item [(ii)] the loop $\psi$ extends to a loop of global
Hamiltonian diffeomorphisms of $W$, contractible in the class of loops
fixing $p$,

\item [(iii)] the loop $\psi$ is contractible in the group of germs of
Hamiltonian diffeomorphisms at $p$,

\item [(iv)] $\mu(\psi)=0$.

\end{itemize}
\end{Lemma}

\begin{proof}
The implications (ii)$\Rightarrow$(i) and (iii)$\Rightarrow$(iv) are clear
and (i)$\Rightarrow$(iv) is established above.

To prove that (iv)$\Rightarrow$(iii), we identify a neighborhood of
$p$ in $W$ with a neighborhood of the origin in $\R^{2n}$. Then, as is
easy to see, the loop $\psi^t$ is homotopy equivalent to its
linearization $d\psi^t$, a loop of (germs of) linear maps.  By the
definition of the Maslov index, $d\psi^t$ is contractible in $\Sp(\R^{2n})$
if and only if $\mu(\psi):=\mu(d\psi)=0$.

To complete the proof of the lemma, it remains to show that
(iii)$\Rightarrow$(ii).

To this end, let us first analyze the case where $\psi^t$ is
$C^1$-close to the identity. Fixing a small neighborhood $U$ of $p$,
we identify a neighborhood of the diagonal in $U\times U$ with a
neighborhood of the zero section in $T^*U$. Then the graphs of
$\psi^t$ in $U\times U$ turn into Lagrangian sections of $T^*U$.
These sections are the graphs of exact forms $df_t$ on $U$, where all
$f_t$ are $C^2$-small and $f_0\equiv 0\equiv f_1$. Then we extend (the
germs of) the functions $f_t$ to $C^2$-small functions $\tf_t$ on $W$
such that $\tf_0\equiv 0\equiv \tf_1$. The graphs of $d\tf_t$ in
$T^*W$ form a loop of exact Lagrangian submanifolds which are
$C^1$-close to the zero section. Thus, this loop can be viewed as a
loop of Hamiltonian diffeomorphisms of $W$. It is clear that the
resulting loop is contractible in the class of loops fixing $p$.

To deal with the general case, consider a family $\psi_s$, $s\in
[0,\,1]$, of local loops with $\psi_0\equiv \id$ and $\psi_1^t
=\psi^t$.  Let $0=s_0<s_1<\cdots<s_k=1$ be a partition of the interval
$[0,\,1]$ such that the loops $\psi_{s_i}$ and $\psi_{s_{i+1}}$ are
$C^1$-close for all $i=0,\ldots, k-1$. In particular, the loop
$\psi_{s_1}$ is $C^1$-close to $\psi_0=\id$, and thus extends to a
contractible loop $\tpsi_{s_1}$ on $W$. Arguing inductively, assume
that a contractible extension $\tpsi_{s_i}$ of $\psi_{s_i}$ has been
constructed.  Consider the loop
$\eta^t=\psi_{s_{i+1}}^t(\psi_{s_i}^t)^{-1}$ defined near $p$. This
loop is $C^1$-close to the identity, for $\psi_{s_{i+1}}$ and
$\psi_{s_i}$ are $C^1$-close. Hence, $\eta$ extends to a contractible
loop $\tilde{\eta}$ on $W$. Then
$\tpsi_{s_{i+1}}^t:=\tilde{\eta}^t\tpsi^t_{s_i}$ is the required
extension of $\psi_{s_{i+1}}$, contractible in the class of loops
fixing $p$.
\end{proof}

\begin{Remark}
\labell{rmk:loop-ext}

It is clear from the proof of Lemma
\ref{lemma:loop-ext} that the extension of the germ of a loop near $p$
to a global loop fixing $p$ can be carried out with some degree of
control of the $C^k$-norm and the support of the loop. We will need
the following simple fact, which can be easily verified by adapting
the proof of the implication (iii)$\Rightarrow$(ii).

Assume that $\psi$ is the germ of a loop near $p$ and
the linearization of $\psi$ at $p$ is equal to the identity:
$d\psi^t_p=I$ for all $t$. Then $\psi$ extends to a loop  $\tpsi$ of global
Hamiltonian diffeomorphisms of $W$ such that $\tpsi$ is contractible in
the class of loops fixing $p$ and having identity linearization at $p$.
\end{Remark}

\section{Local Floer homology}
\labell{sec:LFH}

\subsection{Local Morse homology}
\labell{sec:LMH} Let $f\colon M^m\to\R$ be a smooth function on a
manifold $M$ and let $p\in M$ be an isolated critical point of
$f$. Fix a small neighborhood $U$ of $p$ containing no other critical
points of $f$ and consider a small generic perturbation $\tf$ of $f$
in $U$. To be more precise, $\tf$ is Morse inside $U$ and
$C^1$-close to $f$. Then, as is easy to see, for any
two critical points of $\tf$ in $U$, all anti-gradient trajectories
connecting these two points are contained in $U$. Moreover, the same
is true for broken trajectories connecting these two points. As a consequence,
the vector space (over $\Z_2$) generated by the critical points of
$\tf$ in $U$ is a complex with (Morse) differential defined in the
standard way. (See, e.g., \cite{schwarz:book}.) Furthermore, the
continuation argument shows that the homology of this
complex, denoted here by $\HML(f,p)$ and referred to as the \emph{local
Morse homology} of $f$ at $p$, is independent of the choice of
$\tf$. This construction is a particular case of the one
from \cite{F:witten}.

\begin{Example}
Assume that $p$ is a non-degenerate critical point of $f$ of index $k$.
Then $\HML(f,p)=\Z_2$ when $*=k$ and $\HML(f,p)=0$ otherwise.
\end{Example}

\begin{Example}
\labell{exam:Morse-max}
When $p$ is a strict local maximum of $f$, we have
$\HMLM(f,p)=\Z_2$. Indeed, in this case, as is easy to see from 
standard Morse theory,
$$
\HMLM(f,p)= \Ho_m(\{f\geq f(p)-\eps\}, \{f= f(p)-\eps\})=\Z_2,
$$
where $\eps>0$ is assumed to be small and such that $f(p)-\eps$ is a
regular value of $f$.
\end{Example}

We will need the following two properties of local Morse homology:

\begin{enumerate}
\item[(LM1)] Let $f_s$, $s\in [0,\, 1]$, be a family of smooth
  functions with \emph{uniformly isolated} critical point $p$, i.e.,
  $p$ is the only critical point of $f_s$, for all $s$, in some
  neighborhood of $p$ independent of $s$.  Then $\HML(f_s,p)$ is
  constant throughout the family, and hence $\HML(f_0,p)=\HML(f_1,p)$;
  cf.\ \cite[Lemma 4]{GM}.

\item[(LM2)] The function $f$ has a (strict) local maximum at $p$ if and
only if $\HMLM(f,p)\neq 0$, where $m=\dim M$.
\end{enumerate}

The first assertion, (LM1), is again established by the continuation
argument; cf.\ \cite{schwarz:book}. We emphasize that here the
assumption that $p$ is uniformly isolated is essential and cannot be
replaced by the weaker condition that $p$ is just an isolated critical
point of $f_s$ for all $s$. (Example: $f_s(x)=s x^2+ (1-s)x^3$ on $\R$
with $p=0$. The author is grateful to Doris Hein for this remark.)

Regarding (LM2) first note that, by Example \ref{exam:Morse-max},
$\HMLM(f,p)\neq 0$ when $f$ has a strict local maximum at $p$. The
converse requires a proof although the argument is quite standard.

\begin{proof}[Proof of the implication $(\Leftarrow)$ in (LM2)]
Denote by $\psi^t$ the anti-gradient flow of $f$.  Let $B$ be a
closed connected neighborhood of $p$ with piecewise smooth boundary $\p B$
such that whenever $x\in B$ and $\psi^t(x)\in B$ the entire
trajectory segment $\psi^\tau(x)$ with $\tau\in [0,\,t]$ is also in
$B$, and $p$ is the only critical point of $f$ contained in $B$.
We call $B$ a Gromoll--Meyer neighborhood of $p$. It is not hard to
show that $p$ has an (arbitrarily small) Gromoll--Meyer
neighborhood; see \cite[pp.\ 49--50]{Ch} or \cite{GM}. (Strictly
speaking, the above definition is slightly different from the one
used in \cite{Ch}. However, the existence proof given in \cite[pp.\
49--50]{Ch} goes through with no modifications.)
When $\tf$ is a $C^2$-small generic perturbation of $f$ supported in
$B$, the Morse complex of $\tf|_B$ is defined and its homology is
equal to $\HML(f,p)$. (The fact that $\p^2=0$ follows from the
requirements on $B$.) Assume that $\HMLM(f,p)\neq 0$. Then there exists a
non-zero cycle $C$ of degree $m$ in the Morse complex of $\tf|_B$.
Let $V$ be the closure of the union of the unstable manifolds of
$\tf|_B$ for all local maxima entering $C$. The set $V$ is the closure of 
a domain
with piecewise smooth boundary. The condition that $C$ is a cycle
implies that for every critical point $x$ of $\tf$ in $V$, the intersection
of the unstable manifold of $x$ with $B$ is contained entirely in the
interior of $V$. Hence, $\p V\subset \p B$, and thus $B=V$. It follows that
at every smooth point $z\in \p B$, the gradient $\nabla \tf (z)=\nabla f(z)$
either points inward or is tangent to~$\p B$.

Consider a Gromoll--Meyer neighborhood $N$ of $p$. Note that for a
small generic $\eps>0$ the connected component $B$ of $N\cap \{
f(p)-\eps \leq f \leq f(p)+\eps\}$ containing $p$ is also a
Gromoll--Meyer neighborhood. Clearly, when $p$ is not a local
maximum of $f$, there are smooth points on $\p B$ where $\nabla f$
points inward, provided that $\eps>0$ is small. As a consequence of
the above analysis, $\HMLM(f,p)=0$ if $p$ is not a local maximum.
This completes the proof of the implication $(\Leftarrow)$.
\end{proof}

\begin{Remark}
  Generalizing Example \ref{exam:Morse-max} and the proof of (LM2), it
  is not hard to relate local Morse homology to local homology of a
  function, introduced in \cite{Mo}; see also \cite{Ch,GM}.  However,
  we do not touch upon this question, for such a generalization is not
  necessary for the proof of Theorem \ref{thm:main}.  In the setting
  of local homology, the analogues of (LM1) and (LM2) are established
  in \cite{Ch,GM} and, respectively, in \cite{Hi:grows,Hi}.  
\end{Remark}

\subsection{Local Floer homology: the definition and basic properties}
\labell{sec:LFH2}
Let $\gamma$ be an isolated one-periodic orbit of a
Hamiltonian $H\colon S^1\times W\to \R$. Pick a sufficiently small
tubular neighborhood $U$ of $\gamma$ and consider a non-degenerate
$C^2$-small perturbation $\tH$ of $H$ supported in $U$. More
specifically, let $U$ be a neighborhood of $\gamma(S^1)$, where
$\gamma$ is viewed as a curve in the extended phase space $S^1\times
W$, and let $\tH$ be a Hamiltonian $C^2$-close to $H$, equal to $H$
outside of $U$, and such that all one-periodic orbits of $\tH$ that
enter $U$ are non-degenerate. (Such perturbations $\tH$ do exist; see
\cite[Theorem 9.1]{SZ}.) Abusing notation, we will treat
$U$ simultaneously as an open set in $W$ and in $S^1\times W$.

Consider one-periodic orbits of $\tH$ contained in $U$. Every
anti-gradient trajectory $u$ connecting two such orbits is also contained
in $U$, provided that $\|\tH-H\|_{C^2}$ and $\supp(\tH-H)$ are small
enough. Indeed, the energy $E(u)$ is equal to the difference of action
values on the periodic orbits, and thus is bounded from above by
$O(\|\tH-H\|_{C^2})$. The $C^2$-norm of $\tH$ is bounded from above by
a constant independent of $\tH$, say $2\|H\|_{C^2}$.  Therefore,
$|\p_s u|$ is pointwise uniformly bounded  by
$O(\|\tH-H\|_{C^2})$, and it follows that $u$ takes values in $U$; see
\cite{Sa:london, Sa}.  Note also that for a suitable small
perturbation of a fixed almost complex structure on $W$ the
transversality requirements are satisfied for moduli spaces of Floer
anti-gradient trajectories connecting one-periodic orbits $\tH$
contained in $U$, \cite{FHS,SZ}.

By the compactness theorem, every broken anti-gradient trajectory $u$
connecting two one-periodic orbits in $U$ lies entirely in
$U$. Hence, the vector space (over $\Z_2$) generated by one-periodic
orbits of $\tH$ in $U$ is a complex with (Floer) differential defined
in the standard way. The continuation argument (see, e.g., \cite{SZ})
shows that the homology of this complex is independent of the choice
of $\tH$ and of the almost complex
structure. We refer to the resulting homology group $\HFL(H,\gamma)$
as the \emph{local Floer homology} of $H$ at $\gamma$.
Homology groups of this type were first considered (in a more general
setting) by Floer in \cite{F:witten,Fl}; see also \cite[Section
3.3.4]{Poz}. In fact, an orbit $\gamma$ can be replaced by a connected
isolated set $\Gamma$ of one-periodic orbits of $H$; \cite{F:witten,Fl,Poz}.
(Note that $A_H|_{\Gamma}$ is constant, for $A_H$
is continuous and $\CS(H)$ is nowhere dense.)

\begin{Example}
Assume that $\gamma$ is non-degenerate and $\MUCZ(\gamma)=k$.
Then $\HFL(H,\gamma)=\Z_2$ when $*=k$ and $\HFL(H,\gamma)=0$ otherwise.
\end{Example}

We will need the following properties of local Floer homology:

\begin{enumerate}
\item[(LF1)] Let $H^s$, $s\in [0,\, 1]$, be a family of Hamiltonians
  such that $\gamma$ is a \emph{uniformly isolated} one-periodic orbit
  for $H^s$, i.e., $\gamma$ is the only periodic orbit of $H^s$,
  for all $s$, in some open set independent of $s$.  Then
  $\HFL(H^s,\gamma)$ is constant throughout the family, and hence
  $\HFL(H^0,\gamma)=\HFL(H^1,\gamma)$.
\end{enumerate}

This is again an immediate consequence of the continuation
argument. However, it is worth pointing out that
unless $H^s$ is monotone decreasing, the isomorphism
$\HFL(H^0,\gamma)=\HFL(H^1,\gamma)$ is not induced by the homotopy
$H^s$ in the same sense as the homomorphism $\Psi_{H^0,H^1}$
is induced by a monotone homotopy; see \eqref{eq:homo}.
The isomorphism in question is constructed similarly to \eqref{eq:viterbo}
by breaking $H^s$ into a composition of homotopies close to the identity.

Local Floer homology spaces are building blocks for filtered Floer
homology. Namely, essentially by definition, we have the following

\begin{enumerate}
\item[(LF2)] Let $c\in \R$ be such that all one-periodic orbits $\gamma_i$ of
$H$ with action $c$ are isolated. (As a consequence, there are only finitely
many such orbits.) Then, if $\eps>0$ is small enough,
$$
\HF_*^{(c-\eps,\,c+\eps)}(H)=\bigoplus_i \HFL(H,\gamma_i).
$$
In particular, if all one-periodic orbits $\gamma$ of $H$ are isolated
and ${\operatorname{HF}}^{\mathit{loc}}_{k}(H,\gamma)=0$ for some $k$
and all $\gamma$, we have ${\operatorname{HF}}_{k}(H)=0$ by the
long exact sequence \eqref{eq:seq} of filtered Floer homology.
\end{enumerate}

The effect on local Floer homology of the composition of $H$ with a
loop of Hamiltonian diffeomorphisms is the same as in the global setting
and is established in a similar fashion; see Section
\ref{sec:loops}.

\begin{enumerate}
\item[(LF3)] Let $\psi^t=\varphi^t_G$ be a loop of Hamiltonian
diffeomorphisms of $W$. Then
$$
\HFL(G\# H,\psi(\gamma)) = \HFL(H,\gamma)
$$
for every isolated one-periodic orbit $\gamma$ of $H$, where $\psi(\gamma)$
stands for the one-periodic orbit $\psi^t(\gamma(t))$ of $G\# H$ corresponding
to $\gamma$; see Section \ref{sec:loops}.
\end{enumerate}

As is clear from the definition of local Floer homology, $H$  need not be
a function on the entire manifold $W$ -- it
is sufficient to consider Hamiltonians defined only on a neighborhood of
$\gamma$. For the sake of simplicity, we focus on the particular
case, relevant here, where $\gamma(t)\equiv p$ is a constant orbit,
and hence $dH_t(p)=0$ for all $t\in S^1$. Then (LF1) still holds and
(LF3) takes the following form:

\begin{enumerate}
\item[(LF4)] Let $\psi^t=\varphi^t_G$ be a loop of Hamiltonian diffeomorphisms
defined on a neighborhood of $p$ and fixing $p$ (i.e., $\psi^t(p)=p$ for all
$t\in S^1$). Then
$$
\HFL(G\# H,p) = {\operatorname{HF}}^{\mathit{loc}}_{\ast-2\mu}(H,p),
$$
where $\mu$ is the Maslov index of the loop
$t\mapsto d\psi^t_p\in \Sp(T_pW)$.
\end{enumerate}

Note that in (LF3), in contrast with (LF4), we \emph{a priori} know
that $\mu=0$ as is pointed out in Section \ref{sec:loops}. Hence, the
shift of degrees does not occur when $\psi^t$ is a global loop.  In
other words, comparing (LF3) and (LF4), we can say that the group
$\HFL(H,\gamma)$ is completely determined by the Hamiltonian
diffeomorphism $\varphi_H\colon W\to W$ and its fixed point
$\gamma(0)$, while the germ of $\varphi_H$ at $p$ determines $\HFL(H,p)$
only up to a shift in degree. The degree depends on the class
of $\varphi_H^t$ in the universal covering of the group of germs of
Hamiltonian diffeomorphisms.

Finally note that in the construction of local Floer homology the
Hamiltonian $H$ need not have period one. The definitions and results
above extend word-for-word to $T$-periodic Hamiltonians and, in
particular, to the $T$-th iteration $H^{(T)}$ of a one-periodic
Hamiltonian $H$ as long as the $T$-periodic orbit in question is
isolated.

\subsection{Local Floer homology via local Morse homology}
\labell{sec:LFH-LMH}

A fundamental property of Floer homology is that $\HF_*(H)$
is equal to the Morse homology of a smooth function on
$W$ (and thus to the homology of $W$). The key to establishing this
fact is identifying
$\HF_*(H)$ with $\HM_{*+n}(H)$, when the
Hamiltonian $H$ is autonomous and $C^2$-small; see \cite{FHS,SZ}. A
similar identification holds for local Floer homology. We consider
here the case of $T$-periodic Hamiltonians, for this is the
(superficially more general) situation where the results will be
applied in the subsequent sections.

\begin{Example}
\labell{ex:LFH-LMH}
Assume that $p$ is an isolated critical point of an autonomous Hamiltonian
$F$ and
\begin{equation}
\labell{eq:LFH-LMH1}
 T\cdot \|d^2 F_p\|<2\pi.
\end{equation}
Then $\HFL\big(F^{(T)},p\big)=
{\operatorname{HM}}^{\mathit{loc}}_{\ast+n}(F,p)$.  Indeed, when the
condition \eqref{eq:LFH-LMH1} is satisfied, the Hamiltonians $sF$,
$s\in (0,\, 1]$, have no non-trivial $T$-periodic orbits (uniformly)
near $p$. (See \cite[pp. 184--185]{HZ} or the proof of Lemma
\ref{lemma:LFH-LMH} below.) Thus, $p$ is a uniformly isolated
$T$-periodic orbit of $sF$ for $s\in [\eps,\,1]$ when $\eps>0$ is
small, and $\HFL\big(sF^{(T)},p\big)$ is constant throughout this
family by (LF1).  The argument of \cite{FHS,SZ} shows that the Floer
complex of $sF^{(T)}=sT\cdot F$ is equal to the local Morse complex of
$F$ when $s$ is close to zero.
\end{Example}

In what follows, we will need a slightly more general version of this fact,
where the Hamiltonian is ``close'' to a function independent of time.

\begin{Lemma}
\labell{lemma:LFH-LMH} 
Let $F$ be a smooth function and let $K$ be a
$T$-periodic Hamiltonian, both defined on a neighborhood of a point
$p$. Assume that $p$ is  an isolated critical point of $F$, and the following
conditions are satisfied:
\begin{itemize}

\item The inequalities $\| X_{K_t}-X_F\|\leq \eps \| X_F\|$ and
$\|\dot{X}_{K_t} \|\leq \eps \| X_F\|$ hold pointwise near $p$ for all
$t\in S^1_T$. (The dot stands for the derivative with respect to
time.)

\item The Hessians $d^2 (K_t)_p$ and $d^2 F_p$ and the constant $\eps>0$
are sufficiently small. Namely, $\eps<1$ and
\begin{equation}
\label{eq:LFH-LMH2}
T\cdot \left(\eps(1-\eps)^{-1}+\max_{t}\|d^2 (K_t)_p\|+
\|d^2 F_p\|\right)< 2\pi.
\end{equation}
\end{itemize}
Then $p$ is an isolated $T$-periodic orbit of $K$. Furthermore,
\begin{enumerate}
\item[(a)] $\HFL(K^{(T)},p)=
{\operatorname{HM}}^{\mathit{loc}}_{\ast+n}(F,p)$;
\item[(b)]
if $\HFLN(K^{(T)},p)\neq 0$,
the functions $K_t$ for all $t$ and $F$ have a strict local maximum at $p$.
\end{enumerate}
\end{Lemma}

\begin{Remark}
  The requirement of this lemma, asserting that $K$ is in a certain
  sense close to $F$, plays a crucial role in our proof of Theorem
  \ref{thm:main} (cf.\ Lemmas \ref{lemma:local} and
  \ref{lemma:local2}) and in the argument of \cite{Hi}. To the best of
  the author's knowledge, this requirement is originally introduced in
  \cite[Lemma 4]{Hi} as that $K$ is relatively autonomous.  In what
  follows, we will sometimes call $F$ a \emph{reference function} and
  say that the pair $(F,K)$ meets the requirements of Lemma
  \ref{lemma:LFH-LMH}.
\end{Remark}

\begin{proof}
  Since the statement is local, we may assume that
  $p=0\in\R^{2n}=W$. Consider the family of Hamiltonians
  $K^s=(1-s)K+sF$ starting with $K^0=K$ and ending with $K^1=F$. We
  claim that $\gamma\equiv p$ is a uniformly isolated $T$-periodic
  orbit of $K^s$ for $s\in [0,\, 1]$.

We show this by adapting the proof of \cite[Proposition 17, p.\
184]{HZ}. Fix $r>0$ and let $B_r$ be the ball of radius $r$ centered
at $p$. Since $p$ is a constant $T$-periodic orbit of $K^s$, every
$T$-periodic orbit $\gamma$ of $K^s$ with $\gamma(0)$ sufficiently
close to $p$ is contained in $B_r$.  Recall also that
$2\pi \|z\|_{L^2}\leq T\| \dot{z}\|_{L^2}$ for any map
$z\colon S^1_T\to\R^{2n}$ with zero mean. Applying this inequality to
$z=\dot{\gamma}=X_{K^s}(\gamma)$, we obtain
\begin{equation*}
\begin{split}
\frac{2\pi}{T}\cdot \| \dot{\gamma}\|_{L^2} &\leq \| \ddot{\gamma}\|_{L^2}\\
&= \left\|\frac{d}{dt} X_{K^s}(\gamma) \right\|_{L^2}\\
&\leq \left\| \dot{X}_{K^s}(\gamma) \right\|_{L^2}
+ \left\| \nabla^2 K^s (\gamma)\dot{\gamma} \right\|_{L^2}\\
&\leq \eps\left\| X_{F}(\gamma) \right\|_{L^2}
+ \left(\max_{t}\|d^2 (K^s_t)_p\|+O(r)\right)
\left\|\dot{\gamma} \right\|_{L^2}.
\end{split}
\end{equation*}
Furthermore, from the first requirement on $X_F$ and $X_K$, it is easy
to see that
\begin{equation}
\labell{eq:FH}
\|X_F\|\leq (1-(1-s)\eps)^{-1}\| X_{K^s}\|\leq (1-\eps)^{-1}\| X_{K^s}\|
\end{equation}
pointwise. Hence,
\begin{equation*}
\begin{split}
\frac{2\pi}{T} \cdot \| \dot{\gamma}\|_{L^2}
&\leq \eps(1-\eps)^{-1}\left\| X_{K^s}(\gamma) \right\|_{L^2}
+ \left(\max_{t}\|d^2 (K^s_t)_p\|+O(r)\right)
\left\|\dot{\gamma} \right\|_{L^2}\\
&\leq \eps(1-\eps)^{-1}\left\| \dot{\gamma}\right\|_{L^2}
+ \left(\max_{t}\|d^2 (K^s_t)_p\|+O(r)\right)
\left\|\dot{\gamma} \right\|_{L^2}\\
&=\left(\eps(1-\eps)^{-1}+\max_{t}\|d^2 (K^s_t)_p\|+O(r)\right)
\|\dot{\gamma}\|_{L^2}.
\end{split}
\end{equation*}
Once \eqref{eq:LFH-LMH2} holds and $r>0$ is small, we have
$$
\eps(1-\eps)^{-1}+\max_{t}\|d^2 (K^s_t)_p\|+O(r)<\frac{2\pi}{T}.
$$ 
Therefore,
$\dot{\gamma}=0$. In other words, $\gamma$ is a constant loop, and thus a
critical point of $K^s_t$ for $t\in S^1_T$. Then, by \eqref{eq:FH},
$dF(\gamma)=0$.  As a consequence, $\gamma\equiv p$ since $p$ is an
isolated critical point of $F$. This shows that $p$ is a uniformly isolated
$T$-periodic orbit of $K^s$.

By (LF1), the local Floer homology $\HFL\big((K^s)^{(T)},p\big)$ is
constant throughout the family $K^s$, and
$\HFL\big(K^{(T)},p\big)=\HFL\big(F^{(T)},p\big)$. As a consequence of
\eqref{eq:LFH-LMH2}, the condition \eqref{eq:LFH-LMH1} of Example
\ref{ex:LFH-LMH} is satisfied. Applying this example, we conclude that
$\HFL\big(K^{(T)},p\big)=\HFL\big(F^{(T)},p\big)=
{\operatorname{HM}}^{\mathit{loc}}_{\ast+n}(F,p)$.  This proves (a).

By (LM2), $p$ is an isolated local maximum of $F=K^1_t$, and hence, as
is easy to see from the first condition of the lemma, $p$ is a
uniformly isolated critical point of $K^s_t$ for $s\in [0,\,1]$ and
every fixed $t\in S^1_T$.  Now, by (LM1) and (LM2) applied to
$f_s=K^s_t$, all functions $K^s_t$, and, in particular, $K_t=K^0_t$,
have a (strict) local maximum at $p$.  This proves (b) and concludes the
proof of the lemma.
\end{proof}

\section{Proof of Theorem \ref{thm:main}}
\labell{sec:proof}

As has been pointed out above, it is sufficient to prove the theorem
for (contractible!)  periodic orbits of a Hamiltonian $H$ generating
$\varphi$ rather than for all periodic points of $\varphi$.  Let
$H\colon S^1\times W\to \R$ be a one-periodic Hamiltonian with
finitely many one-periodic orbits $\alpha$. Then, these orbits are
isolated and the action spectrum of $H$ is comprised of finitely many
points.

For every one-periodic orbit $\alpha$ of $H$ denote by
$d_1(\alpha),\ldots,d_{m_\alpha}(\alpha)$ the degrees of roots of
unity, different from 1, among the Floquet multipliers of $\alpha$.

Arguing by contradiction, assume that for every sufficiently large
integer $\tau$, all $\tau$-periodic orbits of $H$ are iterated or, in
other words, $\varphi_H$ has only finitely many \emph{simple periods},
i.e., periods of simple, non-iterated, orbits. In particular, every
periodic orbit of $H$ with sufficiently large period is iterated. Let
$m_1,\ldots,m_k$ be the finite collection of integers comprised of all
simple periods (greater than 1) and the degrees $d_j(\alpha)$ for
all one-periodic orbits $\alpha$. Then, in particular, every
$\tau$-periodic orbit is an iterated one-periodic orbit when $\tau$
which is not divisible by any of the integers $m_j$. Moreover, all
$\tau$-periodic orbits are isolated and $\CS(H^{(\tau)})=\tau\CS(H)$.

Recall also that, when $\alpha$ is a weakly non-degenerate one-periodic orbit 
of $H$ and $\tau$ is a sufficiently large integer, 
not divisible by $d_1(\alpha),\ldots,d_{m_\alpha}(\alpha)$, we have
\begin{equation}
\labell{eq:hom-vanishing}
\HFLN\big(H^{(\tau)},\alpha^{(\tau)}\big)=0.
\end{equation}
Indeed, as is shown in \cite{SZ}, for a generic perturbation of
$H$ supported near $\alpha$, the orbit $\alpha^{(\tau)}$ splits into
non-degenerate orbits with Conley--Zehnder index different from $n$.

Next observe that there exists a strongly degenerate one-periodic
orbit $\gamma$ of $H$ such that $\gamma^{(\tau_i)}$ is an isolated
$\tau_i$-periodic orbit for some sequence $\tau_i\to\infty$ and
\begin{equation}
\labell{eq:hom-nonvanishing}
\HFLN\big(H^{(\tau_i)},\gamma^{(\tau_i)}\big)\neq 0,
\end{equation}
where all $\tau_i$ are divisible by $\tau_1$ and none of $\tau_i$ is
divisible by $m_1,\ldots,m_k$.

To prove this, first note that by \eqref{eq:hom-vanishing} for any
sufficiently large integer $\tau$, not divisible by $m_1,\ldots,m_k$,
there exists a totally degenerate one-periodic orbit $\eta$ such that
$\HFLN\big(H^{(\tau)},\eta^{(\tau)}\big)\neq 0$. (Otherwise,
\eqref{eq:hom-vanishing} held for all $\tau$-periodic orbits, and we
would have $\HF_n\big(H^{(\tau)}\big)=0$ by (LF2).) Pick an infinite
sequence $\tau'_1<\tau'_2<\ldots$ of such integers satisfying the
additional requirement that $\tau'_{i+1}$ is divisible by $\tau'_i$
for all $i\geq 1$. (For instance, we can take $\tau'_i=q^i$, where $q$
is a sufficiently large prime.) As we have observed, for every
$\tau'_i$ there exists a totally degenerate one-periodic orbit
$\eta_i$ such that
$\HFLN\big(H^{(\tau'_i)},\eta_i^{(\tau'_i)}\big)\neq 0$. Since there
are only finitely many distinct one-periodic orbits, one of the orbits
$\gamma$ among the orbits $\eta_i$ and some infinite subsequence
$\tau_i$ in $\tau'_i$ satisfy \eqref{eq:hom-nonvanishing}. (We also
re-index the subsequence $\tau_i$ to make it begin with~$\tau_1$.)

Let $a=A_H(\gamma)$. We will use the orbit $\gamma$ and the sequence $\tau_i$ to prove

\begin{itemize}

\item[\textbf{Claim:}] \emph{For every $\eps>0$ there exists
$T_0$ such that for any $T>T_0$ and some $\delta_T$ in the range
$(0,\,\eps)$, depending on $T$, we have
\begin{equation}
\labell{eq:range1}
\HF^{(T\tau_1a+\delta_T,\,T\tau_1a+\eps)}_{n+1}\big(H^{(T\tau_1)}\big)\neq 0.
\end{equation}
}
\end{itemize}

The theorem readily follows from the claim. Indeed, set
$\CS(H)=\{c_1,\ldots,c_s\}$. Then, if $T>T_0$ is such that $T\tau_1$
is not divisible by $m_1,\ldots,m_k$, we have
$\CS\big(H^{(T\tau_1)}\big)=\{{T\tau_1}c_1,\ldots,{T\tau_1}c_s\}$.
Thus, for any fixed $\eps>0$ and $0<\delta_T<\eps$, the interval
$(T\tau_1a+\delta_T,\,T\tau_1a+\eps)$ contains no action values of
$H^{(T\tau_1)}$ when $T$ is sufficiently large. This contradicts the
claim. (Note that we have used the assumption that $\varphi$ has
finitely many simple periods twice: the first time to find the orbit
$\gamma$ and the sequence $\tau_i$ and the second time to arrive at the
contradiction with the claim.)

To establish the claim, it is convenient to adopt the following

\begin{Definition}
\labell{def:sympl-deg}

A one-periodic orbit $\gamma$ of a one-periodic Hamiltonian $H$ is
said to be a \emph{symplectically degenerate} maximum if there exists
a sequence of loops $\eta_i$ of Hamiltonian diffeomorphisms such that
$\gamma(t)=\eta_i^t(p)$, i.e., $\eta_i$ sends $p$ to $\gamma$, for
some point $p\in W$ and all $i$ and $t$, and such that the Hamiltonians
$K^i$ given by
$$
\varphi^t_H=\eta_i^t\circ\varphi^t_{K^i}
$$
and the loops $\eta_i$ have the following properties:
\begin{itemize}
\item[(K1)] the point $p$ is a strict local maximum of $K_t^i$ for all
$t\in S^1$ and all $i$,
\item[(K2)]
there exist symplectic bases $\Xi^i$ in $T_pW$ such that
$$
\|d^2 (K_t^i)_p\|_{\Xi^i}\to 0\text{ uniformly in $t\in S^1$, and}
$$
\item[(K3)] the linearization of the loop $\eta_i^{-1}\circ \eta_j$
at $p$ is the identity map for all $i$ and $j$ (i.e.,
$d\big((\eta_i^t)^{-1}\circ \eta_j^t\big)_p=I$ for all $t\in S^1$)
and, moreover, the loop $(\eta_i^t)^{-1}\circ \eta_j^t$ is
contractible to $\id$ in the class of loops fixing $p$ and
having the identity linearization at $p$.
\end{itemize}
\end{Definition}

\begin{Remark}
\labell{rmk:sympl-deg}
Regarding (K1) and (K3) note that since $\gamma(t)=\eta_i^t(p)$, the
point $p$ is a fixed point of the flow
$$
\varphi_{K^i}^t= \big(\eta_i^t\big)^{-1}\circ\varphi_{H}^t
$$
of $K^i$, and thus a critical
point of $K^i_t$ for all $t$.  Furthermore, $p$ is also a fixed
point of the loop $\eta_i^{-1}\circ \eta_j$ for all $i$ and $j$, for
$\eta_i^t=\varphi^t_{H}\circ(\varphi^t_{K^i})^{-1}$,
and hence
$(\eta_i^t)^{-1}\circ \eta_j^t=\varphi^t_{K^i}\circ
(\varphi^t_{K^j})^{-1}$.  We refer the reader to Section
\ref{sec:norms} for the definition and discussion of the norm with
respect to a basis, used in (K2).

The Hamiltonians $K^i$ and $H$ have the same time-one flow and there
is a natural one-to-one correspondence between (contractible)
one-periodic orbits of the Hamiltonians. The Hamiltonians $K^i$ can be
chosen so that $K_t^i(p)$ is constant and equal to $c=A_H(\gamma)$. In
what follows, we will always assume that $K^i$ is normalized in this
way. Then the corresponding orbits of $K^i$ and $H$ have equal actions
and, in particular, all Hamiltonians $K^i$ have the same action
spectrum and action filtration; see
Section~\ref{sec:loops}. Symplectically degenerate maxima are further
investigated in \cite{GG:gap}. In particular, it is shown there that
condition (K3) is superfluous; see \cite[Remark 5.5]{GG:gap}. This
fact is not used in the present paper.
\end{Remark}

\begin{Example}
Assume that $H_t$ has a strict local maximum at $p$ (and
$H_t(p)\equiv\const$) and $d^2 (H_t)_p=0$ for all $t$. Then $p$ is a
symplectically degenerate maximum of $H$. Indeed, we can take $K^i=H$
and $\eta_i=\id$ and any fixed symplectic basis as $\Xi^i$.
More generally, vanishing of the Hessian may be
replaced by the condition that $\|d^2 (H_t)_p\|_{\Xi}$ can be made
arbitrarily small by a suitable choice of $\Xi$, cf.\ \cite{Hi}. This
condition is satisfied, for instance, when $H$ is autonomous and all
eigenvalues of the linearization of $X_H$ at $p$ are equal to zero; see Lemma
\ref{lemma:la}.
\end{Example}

\begin{Example}
\labell{ex:sympl-deg}
Assume that $\gamma$ is a symplectically degenerate maximum of $H$. Let
$\tH$ be a Hamiltonian generating the flow $\psi^t\circ\varphi^t_H$, where
$\psi$ is a loop of Hamiltonian diffeomorphisms. Then,
the periodic orbit $\psi(\gamma)(t):=\psi^t(\gamma(t))$ of $\tH$ is
a symplectically degenerate maximum of $\tH$ as is easy to verify.
(In other words,  symplectic degeneracy is a property of
the fixed point $\gamma(0)$ of the time-one map $\varphi_H$.)

For instance, in the notation of Definition \ref{def:sympl-deg}, the
constant orbit $p$ is a symplectically degenerate maximum of each Hamiltonian
$K^i$.
\end{Example}

Now we are in a position to state the two results that we need to complete
the proof of Theorem \ref{thm:main}.  The first result gives a
Floer homological criterion for an isolated, strongly degenerate
orbit $\gamma$ to be a symplectically degenerate maximum, and thus
translates local Floer homological properties of $\gamma$ to
geometrical features of a constant orbit $p$ of Hamiltonians $K^i$.
The second one asserts non-vanishing of the filtered Floer homology
of an iterated Hamiltonian $H^{(T)}$ for an interval of actions
just above the action $T\cdot A_H(\gamma)$, provided that $\gamma$
is a symplectically non-degenerate maximum of $H$. When applied to
the Hamiltonian $H^{(\tau_1)}$ in place of $H$, where $\tau_1$ is as
in the claim, these results will yield the claim.

\begin{Proposition}
\labell{prop1}
Let $\gamma$ be a strongly degenerate isolated one-periodic orbit of
$H$ such that its $l$-th iteration $\gamma^{(l)}$ is also isolated and
\begin{equation}
\labell{eq:hom-max}
\HFLN(H,\gamma)\neq 0\text{ and }
\HFLN\big(H^{(l)},\gamma^{(l)}\big)\neq 0 \text{ for some $l\geq n+1$.}
\end{equation}
Then $\gamma$ is a symplectically degenerate maximum of $H$.
\end{Proposition}

\begin{Remark}
\labell{rmk:comp}
Note that, similarly to Definition \ref{def:sympl-deg}, requirement
\eqref{eq:hom-max} is a condition on the fixed point $\gamma(0)$ of
$\varphi_H$, independent of a particular choice of $H$.
\end{Remark}

\begin{Proposition}
\labell{prop2}
Let $\gamma$ be a symplectically degenerate maximum of $H$ and let
$c=A_H(\gamma)$.
Then for every $\eps>0$ there exists $T_0$ such that
$$
\HF^{(Tc+\delta_T,\,Tc+\eps)}_{n+1}\big(H^{(T)}\big)\neq 0\text{ for
all $T>T_0$ and some $\delta_T$ with $0<\delta_T<\eps$.}
$$
\end{Proposition}

Combining the propositions, we conclude that whenever a strongly
degenerate one-periodic orbit $\gamma$ of $H$ satisfies
the hypotheses of Proposition \ref{prop1},
for every $\eps>0$ there exists $T_0$ such that
\begin{equation}
\labell{eq:range2}
\HF^{(Tc+\delta_T,\,Tc+\eps)}_{n+1}\big(H^{(T)}\big)\neq 0 \text{ for
all $T>T_0$},
\end{equation}
where $c=A_H(\gamma)$ and $0<\delta_T<\eps$.

To prove the claim, first note that although Propositions \ref{prop1}
and \ref{prop2} are stated for one-periodic Hamiltonians,
similar results hold, of course, for
Hamiltonians and orbits of any period.  Thus, consider the
Hamiltonian $H^{(\tau_1)}$ in place of $H$ and the isolated orbit
$\gamma^{(\tau_1)}$ in place of $\gamma$ in \eqref{eq:hom-max} and
\eqref{eq:range2}. Then the requirement \eqref{eq:hom-max} is met due
to \eqref{eq:hom-nonvanishing}: $l=\tau_i/\tau_1\geq n+1$ if $i$ is
large enough, since $\tau_i\to\infty$. Furthermore,
$c=A_{H^{(\tau_1)}}(\gamma^{(\tau_1)})=\tau_1 a$ and \eqref{eq:range1}
follows immediate from \eqref{eq:range2}.

It remains to establish Propositions \ref{prop1} and \ref{prop2} to
complete the proof of the theorem. 

\begin{Remark} It is illuminating to compare the above proof with the
  argument due to Salamon and Zehnder from \cite{SZ} asserting that
  ever large prime is a simple period whenever all one-periodic orbits
  of $H$ are weakly non-degenerate. (In particular, the number of
  simple periods less than or equal to $k$ is of order at least
  $k/\log k$.)  In the context of the present paper relying, of
  course, on \cite{SZ}, this is an immediate consequence of
  \eqref{eq:hom-vanishing}. To be more specific, if $\tau$ is a large
  prime and all $\tau$-periodic orbits are iterated,
  \eqref{eq:hom-vanishing} holds for all weakly non-degenerate
  one-periodic orbits and $\HF_n\big(H^{(\tau)}\big)=0$ by (LF2), if
  there are no totally degenerate one-periodic orbits.  When such
  one-periodic orbits exist, we can no longer use the
  Salamon--Zehnder argument to conclude that every large prime is a
  simple period or even to establish the existence of infinitely many simple
  periods. The reason is that in this case the argument implies that
  for every large prime $\tau$ there is a one-periodic orbit $\gamma$
  such that $\HFLN\big(H^{(\tau)},\gamma^{(\tau)}\big)\neq 0$. It is
  unclear, however, if $\HFLN(H,\gamma)\neq 0$, and hence whether or
  not $\gamma$ is a symplectically degenerate maximum.
\end{Remark}

\section{Proof of Proposition \ref{prop1}}
\labell{sec:pr-prop1}

Our goal in this section is to construct the Hamiltonians $K^i$ and
the loops $\eta_i$ meeting requirements (K1--K3).
This construction relies on two technical lemmas, proved in
Section \ref{sec:gf}, and is carried out in several steps.

First, in Section \ref{sec:reduction}, we reduce the problem to the case where
$\gamma$ is a fixed point $p$ of the flow $\varphi^t_H$.

In Section \ref{sec:near-p}, we construct the Hamiltonians $K^i$ and the
loops $\eta_i$ near $p$. We begin by proving in Section \ref{sec:step2} that
the time one map $\varphi=\varphi^1_H$ can be made $C^1$-close to $\id$ by
an appropriate choice of a canonical coordinate system
$\xi$ near $p$. This is essentially an elementary linear algebra fact
(Lemma \ref{lemma:la}, proved in Section \ref{sec:la}), asserting that
a strongly degenerate linear symplectomorphism can be made arbitrarily
close to the identity by conjugation within the linear symplectic
group.

As a consequence, near $p$, the map $\varphi$ is given by a
generating function $F$ in the coordinate system $\xi$. In Section
\ref{sec:step3}, we show that on a neighborhood of $p$ there exists
a Hamiltonian $K$ with time-one flow $\varphi$, which is in a
certain sense close to $F$. Here, the key result is Lemma
\ref{lemma:local} spelling out the relation between $F$ and $K$ and
established in Section \ref{sec:gf}. Choosing a sequence of coordinate
systems $\xi^i$ so that $\|\varphi-\id\|_{C^1(\xi^i)}\to 0$, we
obtain a sequence of Hamiltonians $K^i$ defined near $p$ and meeting
requirement (K2). Then, again near $p$, the loop $\eta_i$ is defined
by $\eta_i^t=\varphi^t_{H}\circ(\varphi^t_{K^i})^{-1}$.

Utilizing condition \eqref{eq:hom-max}, we show in Section
\ref{sec:step4} that the Maslov index of $\eta_i$ is zero.  This
enables us to relate homological properties of $\gamma\equiv p$ to the
geometrical properties of $K^i$ near $p$ and prove (K1) as a
consequence of Lemma \ref{lemma:LFH-LMH}. (Assertion (K2) easily
follows from the construction of $K^i$.)

Property (K3) is proved in Section \ref{sec:K3}. At this stage, we
further specialize our choice of canonical coordinate systems $\xi^i$
to ensure that all flows $\varphi^t_{K^i}$ have the same linearization
at $p$. Then, the first part of assertion (K3) is obvious. By Lemma
\ref{lemma:loop-ext}, the loops $\eta_i$ extend to $W$, for
$\mu(\eta_i)=0$. This, in turn, gives an extension of $K^i$ to
$W$. Carrying out these extensions with some care, we can guarantee that
(K3) holds in its entirety.

\subsection{Reduction to the case of a constant orbit}
\labell{sec:reduction}
In this section, we reduce the proposition to the case where
\begin{itemize}
\item
\emph{$\gamma\equiv p$ is a constant, strongly degenerate one-periodic
orbit of $H$ and $H_t(p)=c$ for all $t\in S^1$}
\end{itemize}
by constructing a loop of Hamiltonian diffeomorphisms $\psi^t$,
$t\in S^1$, of $W$ such that $\gamma(t)=\psi^t(p)$ with $p=\gamma(0)\in W$.

First recall that for any contractible, closed curve $\gamma\colon
S^1\to W$ there exists a contractible loop of Hamiltonian
symplectomorphisms $\psi^t$ for which $\gamma$ is an integral curve,
i.e., $\gamma(t)=\psi^t(\gamma(0))$; cf.\ \cite[Section 9]{SZ}.

For the sake of completeness, let us outline a proof of this fact.
Consider a smooth family of closed curves
$\gamma_s\colon S^1\to W$, $s\in [0,1]$,
connecting the constant loop $\gamma_0\equiv\gamma(0)$ to $\gamma_1=\gamma$.
It is easy to show that there exists a smooth family of
Hamiltonians $G^s_t$ such that for every $t$, the curve $s\mapsto
\gamma_s(t)$ is an integral curve of $G^s_t$ with respect to $s$,
i.e., $\gamma_s(t)=\varphi_G^s(\gamma(0))$. Let $\psi^t=\varphi^1_G$
be the time-one flow (in $s$) of this family, parametrized by $t\in S^1$. Then
$\gamma(t)=\psi^t(\gamma(0))$.  The family
of Hamiltonians $G^s_t$ can be chosen so that $G_{s,0}\equiv 0\equiv
G_{s,1}$. Then $\psi^t$ is a loop of Hamiltonian diffeomorphisms
with $\psi^0=\id=\psi^1$. As readily follows from the construction, the
loop $\psi$ is contractible.

Composing $\varphi^t_H$ with the loop $(\psi^t)^{-1}$ and adding, if
necessary, a time-dependent constant function to the resulting
Hamiltonian $\hH$, we may assume without loss of generality that
$\gamma(t)\equiv p$ is a fixed point of the flow $\varphi^t_{\hH}
=(\psi^t)^{-1}\varphi^t_H$ for all $t\in S^1$ and $\hH_t(p)\equiv c$.
Then $\hH$ has the same time-one map and the same filtered Floer
homology as $H$. By Example \ref{ex:sympl-deg}
and Remark \ref{rmk:comp}, it is sufficient to prove the proposition
for $\hH$.  Thus, we will assume from now on that $\gamma\equiv p$ and
keep the notation $H$ for the modified Hamiltonian $\hH$.

\subsection{The construction of the Hamiltonians $K^i$ and the loops $\eta_i$
near $p$}
\labell{sec:near-p}

Our main objective in this section is to show that for every
$\sigma>0$, there exists a symplectic basis $\Xi$ in $T_pW$ and a
Hamiltonian $K$ on a neighborhood of $p$ such that the time-one flow
of $K$ is $\varphi$, condition (K1) is satisfied, and
$$
\|d\varphi^t_K|_{T_pW}-I\|_\Xi <\sigma\text{ for all $t\in S^1$.}
$$
Then, clearly, there exists a sequence of symplectic bases $\Xi^i$ and
a sequence of Hamiltonians $K^i$ meeting requirements (K1) and (K2).
The loop $\eta_i$ is defined
near $p$ by $\eta_i^t=\varphi^t_H\circ(\varphi^t_{K_i})^{-1}$.

\subsubsection{Making $\varphi$ close to the identity.}
\labell{sec:step2}
Our first goal is to show that
\begin{itemize}
\item
\emph{for any $\sigma>0$ there exists a symplectic basis $\Xi$ in
$T_pW$ such that $\|d\varphi_p-I\|_\Xi<\sigma$. As a consequence
(cf.\ Example \ref{exam:norm1}), for any $\sigma>0$ there exists a
canonical system of coordinates $\xi$ on a neighborhood $U$ of $p$
such that the $C^1(\xi)$-distance from $\varphi$ to the identity
is less than $\sigma$.}
\end{itemize}

This fact is an immediate consequence of

\begin{Lemma}
\labell{lemma:la} Let $\Phi\colon V\to V$ be a linear symplectic map
of a finite--dimensional symplectic vector space $(V,\omega)$ such
that all eigenvalues of $\Phi$ are equal to one. Then $\Phi$ is
conjugate in $\Sp(V,\omega)$ to a linear map which is arbitrarily
close to the identity.
\end{Lemma}

Indeed, since $p$ is a strongly degenerate fixed point of $H$, all
eigenvalues of $d\varphi_p$ are equal to one. Thus, the desired statement
follows from this lemma applied to $\Phi=d\varphi_p$.  The proof of
the lemma is elementary and provided for the sake of completeness
in Section \ref{sec:la}. Here we only mention that $\Phi$ is given by
an upper triangular matrix in some basis $\Xi$ and, by scaling the
elements of $\Xi$ appropriately, one can make $\Phi$ arbitrarily
close to the identity, cf.\ Example \ref{exam:norm1}.  Hence, we only
need to show that $\Xi$ and the scaling can be made symplectic.

\subsubsection{The Hamiltonian $K$ near $p$}
\labell{sec:step3} Pick a system $\xi$ of canonical coordinates near
$p$ such that $\varphi$ is $C^1(\xi)$-close to the identity. In
particular, $\|d\varphi_p-I\|_{\xi_p}$ is small.  Furthermore, the
map $\varphi$ is given, near $p$, by a generating function $F$. The
precise definition of $F$ and the relation between $F$ and
$\varphi$ and $\xi$ are immaterial at the moment and these issues
will be discussed in Section \ref{sec:gf}. At this stage, we only
need to know that $F$ is defined on a neighborhood of $p$ and
uniquely determined by $\xi$ and $\varphi$. (To make this statement
accurate, let us agree that a canonical coordinate system is
comprised of \emph{ordered} pairs of functions $(x_1,y_1),\ldots,
(x_n,y_n)$ such that $\omega=\sum dx_i\wedge dy_i$. Thus, each
coordinate function is assigned to either $x_i$- or $y_i$-group.)
Moreover, $F$ has the following properties:

\begin{enumerate}
\item[(GF1)] $p$ is an isolated critical point of $F$,
\item[(GF2)] $\|F\|_{C^2(\xi)}=O(\|\varphi-\id\|_{C^1(\xi)})$ and
$\| d^2 F_p\|_{\xi_p}=\|d\varphi_p-I\|_{\xi_p}$.
\end{enumerate}

The second item (GF2) requires, perhaps, a clarification.  First note
that $\|\varphi-\id\|_{C^1(\xi)}$ stands here for the
$C^1(\xi)$-distance from $\varphi$ to $\id$; see Section
\ref{sec:norms}. Furthermore, $F$ and $\|\varphi-\id\|_{C^1(\xi)}$
depend on $\xi$. Thus, in (GF2), we view both $\|F\|_{C^2(\xi)}$ and
$\|\varphi-\id\|_{C^1(\xi)}$ as functions of $\xi$ with $\varphi$
fixed and the second item asserts that $\|F\|_{C^2(\xi)}\leq
\const\cdot\|\varphi-\id\|_{C^1(\xi)}$, where $\const$ is independent of
$\xi$, provided that $\|\varphi-\id\|_{C^1(\xi)}$ is small enough.

More generally, let $f$ and $g$ be non-negative functions of $\xi$ and
some (numerical) variables. We write $f=O(g)$, when $f\leq \const\cdot
g$ pointwise, where $\const$ is independent of $\xi$.  The notation
$f=O_\xi(g)$ will be used when $f\leq \const(\xi)\cdot g$ pointwise
as functions of other variables, with $\const(\xi)$ depending on $\xi$
and possibly becoming arbitrarily large.  Furthermore, we denote by
$B_r(\xi)$ the ball of radius $r$ with respect to $\xi$ centered at
$p$.

We will prove

\begin{Lemma}[\cite{Hi}]
\label{lemma:local}
Let $\xi$ be a coordinate system near $p$ such
that $\|\varphi-\id\|_{C^1(\xi)}$ is small. Then for every sufficiently
small $r>0$ (depending on $\xi$), there exists a one-periodic Hamiltonian
$K_t$ on $B_r(\xi)$ such that
\begin{enumerate}
\item[(i)] the time one-flow $\varphi_K$ of $K$ is $\varphi$,

\item[(ii)] $p$ is an isolated critical point of $K_t$
and $K_t(p)\equiv c$,

\item[(iii)] $\|d^2 (K_t)_p\|_{\xi_p}=O(\|d\varphi_p-I\|_{\xi_p})$,

\item[(iv)] the following estimates  hold pointwise
near $p$:
$$
\|X_K-X_F\|_\xi\leq 
\big(O(\|d^2 F_p\|_{\xi_p})+O_\xi(r)\big)\cdot \|X_F\|_{\xi}
$$
and
$$
\|\dot{X}_K\|_\xi\leq 
\big(O(\|d^2F_p\|_{\xi_p})+O_\xi(r)\big)\cdot \|X_F\|_{\xi},
$$
where the dot denotes the time derivative of a vector field.
\end{enumerate}
\end{Lemma}

Note that in (iv) we could have written $\|d\varphi_p-I\|_{\xi_p}$ in place
of $\|d^2 F_p\|_{\xi_p}$ by (GF2).
The important point here is that $\|d\varphi_p-I\|_{\xi_p}$ and
$\|d^2 F_p\|$ can be made arbitrarily small by choosing an
appropriate coordinate system $\xi$. Then, shrinking the domain of $K$, we
can also make the right hand sides in the estimates (iii) and (iv)
arbitrarily small.

A proof of Lemma \ref{lemma:local} can be extracted from \cite{Hi}.
However, to make our proof of Theorem \ref{thm:main} self-contained,
we provide a detailed argument.  Deferring this to Section
\ref{sec:gf}, we proceed with the proof of Proposition \ref{prop1}.

\subsubsection{Properties (K1) and (K2)}
\labell{sec:step4}
Let $K$ be a Hamiltonian on a neighborhood of $p$, such that (i)--(iv) of Lemma
\ref{lemma:local} are satisfied and $\|\varphi-\id\|_{C^1(\xi)}$ is small.
Our first goal is to prove that $K$ meets requirements (K1) and (K2).

Since $\|d\varphi_p-I\|_{\xi_p}\leq \|\varphi-\id\|_{C^1(\xi)}$, by
(iii), we have
\begin{equation}
\labell{eq:dK2varphi}
\|d^2 (K_t)_p\|_{\xi_p}=O(\|\varphi-\id\|_{C^1(\xi)}),
\end{equation}
and hence (K2) is satisfied when $\|\varphi-\id\|_{C^1(\xi)}$ is
sufficiently small.

To establish (K1), consider the loop
$\eta^t=\varphi_H^t(\varphi_K^t)^{-1}$, where $t\in \R$. Thus,
$$
\varphi^t_H=\eta^t\varphi^t_K\text{ for all $t\in \R$.}
$$
Note that $\eta^1=\id$, i.e., $\eta^t$ with $t\in S^1$ is a loop of
Hamiltonian symplectomorphisms near $p$. We denote this loop by
$\eta|_{S^1}$. The $T$-th iteration $\eta|_{S^1_T}$ of $\eta|_{S^1}$
is simply $\eta^t$ with $t\in S^1_T$.

First, let us prove that the Maslov index $\mu=\mu(\eta|_{S^1})$ of
$\eta|_{S^1}$ is necessarily zero, when $\|\varphi-\id\|_{C^1(\xi)}$ is
small.

Let $\bar{K}$ be the time-average of $K_t$, i.e.,
$$
\bar{K}=\int_0^1 K_t\, dt.
$$
A straightforward calculation utilizing Lemma \ref{lemma:local} and
\eqref{eq:dK2varphi} and detailed in Section \ref{sec:barK-K} shows
that the requirements of Lemma \ref{lemma:LFH-LMH} are met, for any
fixed $T$, by the pair $(\bar{K},K)$, provided that
$\|\varphi-\id\|_{C^1(\xi)}$ is small enough. In other words, these
requirements are satisfied when $\varphi$ is $C^1(\xi)$-close to the
identity and $\bar{K}$ is taken as the reference function in Lemma
\ref{lemma:LFH-LMH} (denoted there by $F$). In particular, $p$ is an
isolated $l$-periodic orbit. We set, $T=l$, where $l$ is as in
\eqref{eq:hom-max}.

Then, by Lemma \ref{lemma:LFH-LMH}(a),
$$
\HFL(K^{(l)},p)={\operatorname{HM}}^{\mathit{loc}}_{\ast+n}(\bar{K},p),
$$
and hence
\begin{equation}
\labell{eq:k2}
{\operatorname{HF}}^{\mathit{loc}}_{k}(K^{(l)},p) = 0
\text{ whenever $|k|>n$.}
\end{equation}

Next note that $\mu(\eta|_{S^1_T})=\mu T$ for any $T\in\Z$.
By (LF4),
$$
{\operatorname{HF}}^{\mathit{loc}}_{n-2\mu T}\big(K^{(T)},p\big)
={\operatorname{HF}}^{\mathit{loc}}_{n}\big(H^{(T)},p\big)\neq 0.
$$
as long as $p$ is an isolated one-periodic orbit of $H^{(T)}$.
Applying this identity to $T=l$, we conclude
from \eqref{eq:hom-max} that
$$
{\operatorname{HF}}^{\mathit{loc}}_{n-2\mu l}\big(K^{(l)},p\big)
={\operatorname{HF}}^{\mathit{loc}}_{n}\big(H^{(l)},p\big)\neq 0.
$$
In particular, since $l\geq n+1$,
\begin{equation}
\labell{eq:k1}
{\operatorname{HF}}^{\mathit{loc}}_{k}\big(K^{(l)},p\big) \neq 0
\text{ for some $k$ with $|k|>n$ if $\mu\neq 0$.}
\end{equation}
Combining \eqref{eq:k2} and \eqref{eq:k1}, we conclude that $\mu=0$. A
different proof of this fact, relying on the properties of the mean
Conley--Zehnder index (see \cite{SZ}), can be found in \cite[Section
5.2]{GG:gaps}.

Recall that the condition $\HFLN(H,p)\neq 0$ is a part of the
assumption \eqref{eq:hom-max} in Proposition \ref{prop1}.  Using (LF4)
again -- this time for $t\in [0,\, 1]$ -- and taking into account that
$\mu=0$, we see that
$$
\HFLN(K,p)=\HFLN(H,p)\neq 0.
$$

Furthermore, when $\varphi$ is sufficiently $C^1(\xi)$-close to the
identity, the requirements of Lemma \ref{lemma:LFH-LMH} with $T=1$ and
the Hamiltonians $K$ and $F$ as in Lemma \ref{lemma:local} are
obviously met due to (GF2), \eqref{eq:dK2varphi}, and (iv). Thus, by
Lemma \ref{lemma:LFH-LMH}(b), the function $K_t$ has strict local
maximum at $p$. This proves (K1).

Applying this construction to a sequence of symplectic bases $\Xi^i$
in $T_pW$ such that $\|d\varphi_p-I\|_{\Xi^i}\to 0$, we obtain a
sequence of Hamiltonians $K^i$, meeting requirements (K1) and (K2),
and also the loops $\eta_i$. We emphasize that $K^i$ and $\eta^i$ have
so far been defined only on a neighborhood of $p$.

\subsubsection{The pair $(\bar{K},K)$} 
\label{sec:barK-K}
The goal of this auxiliary section, which is included for the
sake of completeness, is to show that, as stated above, the pair
$(\bar{K},K)$ satisfies the hypotheses of Lemma \ref{lemma:LFH-LMH}
with $T$ fixed.

To this end, note first that by Lemma \ref{lemma:local}(iv), we have
$$ 
\|X_K-X_F\|\leq \big(O(\|d^2 F_p\|)+O(r)\big)\|X_F\|\text{ and }
\|\dot{X}_K\|\leq \big(O(\|d^2 F_p\|)+O(r)\big)\|X_F\| 
$$ 
pointwise near $p$. (Here the coordinate system $\xi$ is suppressed in
the notation.)  Let us integrate the first of these inequalities with
respect to $t$ over $S^1_T$. Then, since $F$ is independent of time,
we have, again pointwise near $p$,
$$
\|X_{\bar{K}}-X_F\|\leq \int_{S^1} \|X_K-X_F\|\, dt \leq a\|X_F\|
$$
with $a=O(\|d^2 F_p\|)+O(r)$. Thus,
$$
\|X_{\bar{K}}-X_F\|\leq a\|X_F\|
$$
and, as a consequence,
$$
\|X_{\bar{K}}\|\geq \|X_F\|- a\|X_F\|=(1-a)\|X_F\|.
$$
Then 
\begin{eqnarray*}
\|X_K-X_{\bar{K}}\| &\leq& \|X_K-X_F\|+\|X_F-X_{\bar{K}}\|\\
                   &\leq& a\|X_F\|+a\|X_F\|\\
                   &\leq& 2a(1-a)^{-1}\|X_{\bar{K}}\|
                 \end{eqnarray*}
Likewise, 
$$
\|\dot{X}_K\|\leq a(1-a)^{-1} \|X_{\bar{K}}\|.
$$
Therefore,
$$
\|X_K-X_{\bar{K}}\|\leq \eps \|X_{\bar{K}}\| \text{ and
}\|\dot{X}_K\|\leq \eps\|X_{\bar{K}}\|,
$$
where $\eps=2a(1-a)^{-1}$ and all inequalities are pointwise.

Recall now that we can make $a>0$ arbitrarily small (with $T$ fixed)
by making a suitable choice of $\xi$ and then requiring $r>0$ to be
sufficiently small. It follows that we can also make $\eps>0$
arbitrarily small. In the same vein, the left hand side of
\eqref{eq:LFH-LMH2} can be made arbitrarily small.  Furthermore, since
$p$ is an isolated critical point of $F$, it is also an isolated
critical point of $\bar{K}$. Therefore, the pair $(\bar{K},K)$ satisfies 
the hypotheses of Lemma \ref{lemma:LFH-LMH}.

\begin{Remark} 
  As has been pointed out in Section \ref{sec:step4}, a pair of
  function satisfying the hypotheses of Lemma \ref{lemma:local} also
  satisfies the hypotheses of Lemma \ref{lemma:LFH-LMH}. We have shown
  that $(\bar{K},K)$ satisfies the conditions of Lemma
  \ref{lemma:LFH-LMH} whenever $(F,K)$ meets the requirements of Lemma
  \ref{lemma:local}. Moreover, by arguing as in
  this section, it is not hard to show that $(\bar{K},K)$ satisfies the
  conditions of Lemma \ref{lemma:local} (and hence of Lemma
  \ref{lemma:LFH-LMH}) once $(F,K)$ does. We omit this
  (straightforward) calculation, for it is never used in the proof.
\end{Remark}

\subsection{Property (K3) and the extension to $W$}
\labell{sec:K3}
To ensure that (K3) holds, we need to impose an addition requirement on
the bases $\Xi^i$. We will prove

\begin{Lemma}
\labell{lemma:special} There exists a sequence of symplectic bases
$\Xi^i$ in $T_pW$ such that $\|d\varphi_p-I\|_{\Xi^i}\to 0$ and the
flows $\varphi^t_{K^i}$ have the same linearization at $p$.
\end{Lemma}

Here $K^i$ is the sequence of Hamiltonians constructed in Section
\ref{sec:step4} using Lemma \ref{lemma:local}.  We prove Lemma
\ref{lemma:special} in Section \ref{sec:gf} along with
Lemma \ref{lemma:local}.  At this point, we only note that, as will
become clear in Section \ref{sec:gf}, the linearized flow
$d(\varphi^t_{K^i})_p$ is completely determined by $\varphi$ and the
basis $\Xi^i$. In particular, the linearization is independent of the
extension of $\Xi^i$ to a canonical coordinate system $\xi^i$ near
$p$. (Here, we use a convention similar to that of Section
\ref{sec:step3} for canonical coordinate systems: a symplectic basis
is divided into two groups of $n$ vectors spanning Lagrangian
subspaces and this division is a part of the structure of a symplectic
basis.)

Since $\eta^t_i=\varphi_H^t(\varphi_{K_i}^t)^{-1}$, we conclude from
Lemma \ref{lemma:special} that
$$
d\big((\eta_i^t)^{-1}\circ \eta_j^t\big)_p
=d\big(\varphi_{K_i}^t\big)_p^{-1}\circ d\big(\varphi_{K_j}^t\big)_p=I.
$$

Let us now extend the loops $\eta_i$ and the Hamiltonians $K^i$ to
$W$ so that the remaining part of requirement (K3) is met:
the loop $\eta_i^{-1}\circ \eta_j$ is contractible to
$\id$ in the class of loops with identity linearization at $p$.

Recall that the Maslov index of the loop $\eta_i$ is zero, as is shown
in Section \ref{sec:step4}. Hence, by Lemma \ref{lemma:loop-ext}, each
of these loops extends to a loop of Hamiltonian diffeomorphisms of
$W$, contractible in the class of loops fixing $p$. Let us fix such an
extension for $\eta_1$. For the sake of simplicity we denote this
extension by $\eta_1$ again.  Consider now the loop
$\psi_i^t=\big(\eta_1^t\big)^{-1}\eta_i^t$.  Then
$d(\psi^t_i)_p=I$. Hence, by Lemma \ref{lemma:loop-ext} and Remark
\ref{rmk:loop-ext}, $\psi_i$ extends to a loop of Hamiltonian
diffeomorphisms of $W$, contractible in the class of loops with
identity linearization at $p$. Keeping the notation
$\psi_i$ for this extension, we set $\eta_i^t=\eta_1^t\psi_i^t$. It is
clear that $\eta_i$ is contractible in the class of loops with
identity linearization at $p$.

\subsection{Proof of Lemma \ref{lemma:la}}
\labell{sec:la}

Lemma \ref{lemma:la} is an immediate consequence of the following stronger
result which is also used in the proof of Lemma \ref{lemma:special}.

\begin{Lemma}
\labell{lemma:la2} Let $\Phi\colon V\to V$ be a linear symplectic map
of a finite--dimensional symplectic vector space $(V,\omega)$. Assume
that all eigenvalues of $\Phi$ are equal to one. Then $V$ can be
decomposed as a direct sum of two Lagrangian subspaces $L$ and $L'$
with $\Phi(L)=L$. Moreover, by a suitable choice of
$\Psi\in\Sp(V,\omega)$ preserving the subspaces $L$ and $L'$, the map
$\Psi\Phi\Psi^{-1}$ can be made arbitrarily close to the identity.
\end{Lemma}

\begin{proof}
We prove the lemma by induction in $\dim V$.  The statement is obvious
when $V$ is two--dimensional. When $\dim V> 2$, we have the following
alternative:

\begin{itemize}
\item either $K=\ker (\Phi-I)$ contains a symplectic subspace $V_0$
\item or $K=\ker (\Phi-I)$ is isotropic.
\end{itemize}

In the former case, we decompose $V$ as $V_0\oplus V_0^\omega$, where the
superscript $\omega$ denotes the symplectic orthogonal. It is easy to see
that this decomposition is preserved by $\Phi$ and $\Phi|_{V_0}=I_{V_0}$. Now
the assertion follows from the induction hypothesis applied to
$\Phi|_{V_0^\omega}$.

In the latter case, pick a symplectic subspace $V_0$ complementary to
$K=\ker (\Phi-I)$ in $K^\omega$ and an isotropic subspace $N$
complementary to $K^\omega$ in $V$. (We are assuming at the moment
that $V_0\neq \{0\}$, i.e., $L$ is not Lagrangian.) Thus, $V=K^\omega
\oplus N$ and $K^\omega=K\oplus V_0$.  Furthermore, $K$ and $K^\omega$
are preserved by $\Phi$; the spaces $V_0$ and $N$ can be canonically
identified with $K^\omega/K$ and $K^*=V/K^\omega$, respectively; and
$\Phi|_K=I_K$. Note that $\Phi$ induces a symplectic linear map
$\Phi_0\colon V_0\to V_0$ with all eigenvalues equal to one and the
identity map $I_N$ on $N=V/K^\omega$. Hence, using the decomposition
$$
V=K\oplus V_0\oplus N,
$$
we can write $\Phi$ in the block upper-triangular form
$$
\Phi=
\begin{bmatrix}
I_K & A & C\\
0 & \Phi_0 & B\\
0 & 0 & I_N
\end{bmatrix},
$$
where $A\colon V_0\to K$ and $C\colon N\to K$ and $B\colon N\to V_0$.
(There are relations between these operators, resulting from the fact
that $\Phi$ is symplectic.)

Consider a block-diagonal symplectic linear transformation of the
form
$$
\Psi=
\begin{bmatrix}
\Lambda & 0 & 0\\
0 & \Psi_0 & 0\\
0 & 0 & (\Lambda^*)^{-1}
\end{bmatrix},
$$
where $\Psi_0\colon V_0\to V_0$ is symplectic, $\Lambda\colon K\to K$
is invertible, and we have identified $N$ with $K^*$. Then
$$
\Psi\Phi\Psi^{-1}=
\begin{bmatrix}
I_K & \Lambda A \Psi_0^{-1} & \Lambda C \Lambda^* \\
0 & \Psi_0 \Phi_0 \Psi_0^{-1} & \Psi_0 B \Lambda^* \\
0 & 0 & I_N
\end{bmatrix}.
$$
By the induction assumption, there exists a decomposition
$V_0=L_0\oplus L_0'$, where $\Phi_0(L_0)=L_0$, and transformations
$\Psi_0$ preserving this decompositions and making $\Psi_0 \Phi_0
\Psi_0^{-1}$ arbitrarily close to $I_{V_0}$.  Set $L=K\oplus L_0$
and $L=L'_0\oplus N$. Then $\Phi(L)=L$ and the decomposition
$V=L\oplus L'$ is preserved by $\Psi$. Furthermore, noticing that
$\Lambda^*$ is close to zero when $\Lambda$ is close to zero, we can pick
$\Lambda$ to make the off-diagonal entries in $\Phi$ arbitrarily
small. With this choice of $\Psi$, the map $\Psi\Phi\Psi^{-1}$ is
close to $I_V$ if $\Psi_0 \Phi_0 \Psi_0^{-1}$ is close to $I_{V_0}$.

When $K=\ker(\Phi-I)$ is Lagrangian (i.e., $V_0=\{0\}$), no
induction reasoning is needed.  We simply set $L=K$ and let $L'=N$
be an arbitrary complementary Lagrangian subspace.  Then the map
$\Phi$ is decomposed as a two-by-two block upper-triangular matrix,
and, similarly to the argument above, $\Lambda$ is chosen to make
the off-diagonal block arbitrarily small.
\end{proof}

\section{The generating function $F$ and the proofs of Lemmas \ref{lemma:local}
and \ref{lemma:special}}
\labell{sec:gf}

\subsection{Generating functions}
In this section, we recall the definition  of a generating function
on $\R^{2n}$ and set the stage for proving Lemma \ref{lemma:local}.
The material reviewed here is absolutely standard -- it goes back to
Poincar\'e -- and we refer the reader to \cite[Appendix 9]{Ar}
and \cite{We71,We77} for a more detailed discussion of generating functions.

Let us identify $\R^{2n}$ with the Lagrangian diagonal $\Delta\subset
\R^{2n}\times \bar{\R}^{2n}$ via the projection to the first factor,
where $\R^{2n}\times \bar{\R}^{2n}$ is equipped with the symplectic
structure $(\omega,-\omega)$, and fix a Lagrangian complement $N$ to
$\Delta$.  Thus, $\R^{2n}\times \bar{\R}^{2n}$ can now be treated
as~$T^*\Delta$.

Let $\varphi$ be a Hamiltonian diffeomorphism defined on a
neighborhood of the origin $p$ in $\R^{2n}$ and such that
$\|\varphi-\id\|_{C^1}$ is sufficiently small. Then the graph $\Gamma$
of $\varphi$ is close to $\Delta$, and hence $\Gamma$ can be viewed as
the graph in $T^*\Delta$ of an exact form $dF$ near
$p\in\Delta=\R^{2n}$. (We normalize $F$ by $F(p)=0$.) The function
$F$, called the generating function of $\varphi$, has the following
properties:
\begin{enumerate}
\item[(GF1$^\prime$)] $p$ is an isolated critical point of $F$ if and only if
$p$ is an isolated fixed point of $\varphi$,

\item[(GF2$^\prime$)] $\|F\|_{C^2}=O(\|\varphi-\id\|_{C^1})$ and
$\| d^2 F_p\|=\|d\varphi_p-I\|$.
\end{enumerate}
For instance, it is clear that the critical points of $F$ are in one-to-one
correspondence with the fixed points of $\varphi$. If $p$ (the
origin) is an isolated fixed point of $\varphi$, the origin is also an
isolated critical point of $F$. Hence, (GF1$^\prime$) holds. The second
property of $F$, (GF2$^\prime$), is also easy to check;
see the references above.

The function $F$ depends on the choice of the Lagrangian complement
$N$ to $\Delta$. To be specific, we take as $N$ the linear subspace of
vectors of the form $((x,0),(0,y))$ in $\R^{2n}\times \bar{\R}^{2n}$,
where $x=(x_1,\ldots, x_n)$ and $y=(y_1,\ldots,y_n)$ are the standard
canonical coordinates on $\R^{2n}$, i.e., $\omega=\sum dy_i\wedge dx_i$.

In the setting of Section \ref{sec:step3}, let $\xi$ be a coordinate
system near $p\in W$. Using $\xi$, we identify a neighborhood of $p$
in $W$ with a neighborhood of the origin in $\R^{2n}$, keeping the
notation $p$ for the origin.  With this identification, $\varphi$
defined near $p\in W$ turns into a Hamiltonian diffeomorphism
$\xi\varphi\xi^{-1}$ defined near the origin $p\in \R^{2n}$.  By
definition,
$\|\varphi-\id\|_{C^1(\xi)}=\|\xi\varphi\xi^{-1}-\id\|_{C^1}$.
Abusing notation, we denote the resulting Hamiltonian diffeomorphism
$\xi\varphi\xi^{-1}$ near $p\in \R^{2n}$ by $\varphi$ again.  By our
background assumptions, $p$ is an isolated fixed point of $\varphi$,
and thus (GF1) and (GF2) follow immediately from (GF1$^\prime$) and
(GF2$^\prime$), respectively.

Furthermore, Lemma \ref{lemma:local} is an immediate consequence of

\begin{Lemma}[\cite{Hi}]
\labell{lemma:local2} Let $\varphi$ be a Hamiltonian diffeomorphism
of a neighborhood of the origin $p\in\R^{2n}$. Assume that
$p$ is an isolated fixed point of $\varphi$ and
$\|\varphi-\id\|_{C^1}$ is so small that the generating function
$F$ is defined. Then for every sufficiently small $r>0$ (depending on
$\varphi$), there exists a one-periodic Hamiltonian $K_t$ on the ball
$B_r$ of radius $r$ centered at $p$ such that
\begin{enumerate}
\item[(i)] the time one-flow $\varphi_K$ of $K$ is $\varphi$,

\item[(ii)] $p$ is an isolated critical point of $K_t$ and $K_t(p)=0$
for all $t\in S^1$,

\item[(iii)] $\|d^2 (K_t)_p\|=O(\|d^2F_p\|)$,

\item[(iv)] the upper bounds 
\begin{equation}
\labell{eq:KF1}
\|X_K-X_F\|\leq 
\big(O(\|d^2 F_p\|)+O_\varphi(r)\big)\cdot \|X_F\|
\end{equation}
and
\begin{equation}
\labell{eq:KF2}
\|\dot{X}_K\|\leq 
\big(O(\|d^2F_p\|)+O_\varphi(r)\big)\cdot \|X_F\|
\end{equation}
hold pointwise near $p$.
\end{enumerate}
\end{Lemma}

The notation used here is similar to that of Section
\ref{sec:step3}. For instance, \eqref{eq:KF1} should be read as that its left
hand side is pointwise bounded from above by
$\big(C_1\|d^2F_p\|+C_2(\varphi)r\big)\cdot
\|X_F\|$, where $C_1$ is independent of $\varphi$ and $C_2(\varphi)$
depends on $\varphi$ and can be arbitrarily large.

Although the proof of Lemma \ref{lemma:local2} is essentially contained
in \cite{Hi}, for the sake of completeness we give a detailed argument here.

\subsection{Proof of Lemma \ref{lemma:local2}}

The proof of the lemma is organized as follows. First we consider the
time-dependent Hamiltonian $\tK$ generating the flow $\varphi^t$ given
by the family of generating functions $F_t=tF$, $t\in [0,\, 1]$, and
verify (i)--(iv) for $\tK$. The time-one flow of $\tK$ is
$\varphi$. However, in general, the Hamiltonian $\tK$ is \emph{not}
periodic in time. Hence, as the next step, we modify $\tK$ to obtain
the required periodic Hamiltonian $K$ and then again check that the
new Hamiltonian $K$ satisfies (i)--(iv).

\subsubsection{The Hamiltonian $\tK$; properties (i) and (ii)}
Consider the family of generating functions $F_t=tF$ with $t\in
[0,\,1]$.  The family of graphs $\Gamma_t$ of $dF_t$ in $T^*\Delta$
beginning with $\Gamma_0=\Delta$ and ending with $\Gamma_1=\Gamma$
can be viewed as a family of graphs of Hamiltonian diffeomorphisms
$\varphi^t$ near $p$ with $\varphi^0=\id$ and $\varphi^1=\varphi$.
Thus, $\varphi^t$ is a time-dependent Hamiltonian flow with the
time-one map $\varphi$, defined near $p$. Let $\tK_t$ be the
Hamiltonian
generating this flow, normalized by $\tK_t(p)=0$.  

Condition (i) is satisfied for $\tK$ by definition, and (ii) is an
immediate consequence of (GF1$^\prime$) and \eqref{eq:KF1}.  Below, we will
also give a direct proof of (ii).

\subsubsection{The Hamiltonian vector field $X_{\tK}$}
Set $x=(x_1,\ldots,x_n)$ and $y=(y_1,\ldots,y_n)$, and
$$
\p_1 =(\p_{x_1},\ldots,\p_{x_n})\text{ and } \p_2=(\p_{y_1},\ldots,\p_{y_n}).
$$

Let $(\bx^t,\by^t)=\varphi^t(x,y)$. Then,  as is well known
and can be checked by a simple calculation, we have
\begin{equation}
\labell{eq:varphi-F}
\begin{cases}
\bx^t-x &= -\p_2 F_t (\bx^t,y)\\
\by^t-y &= \p_1 F_t (\bx^t,y)
\end{cases}
\end{equation}
Differentiating with respect to time, we obtain the following expression
for the Hamiltonian $X_{\tK}$ (cf.\ \cite{Hi}):
\begin{equation}
\labell{eq:XK1}
X_{\tK_t}(\bx^t,\by^t)=A_t(\bx^t,y)X_F(\bx^t,y),
\end{equation}
where
$$
A_t(x,y)=
\begin{bmatrix}
(I+\p_{12}F_t(x,y))^{-1} & 0\\
(I+\p_{12}F_t(x,y))^{-1}\p_{11}F_t(x,y) & I
\end{bmatrix}.
$$
Here $\p_{12}F_t$ stands for the matrix of partial derivatives
$\p^2F_t/(\p x_i\p y_j)$ and, similarly, $\p_{22}F_t$ is the matrix
$\p^2F_t/(\p y_i\p y_j)$.

In other words, introducing the auxiliary diffeomorphism $\kappa^t$
sending $(\bx^t,\by^t)$ to $(\bx^t,y)$, we can rewrite
\eqref{eq:XK1} as
\begin{equation}
\labell{eq:XK2}
X_{\tK_t}(z)=A_t(\kappa^t(z))X_F(\kappa^t(z))
\end{equation}
for every $z$ near the origin.

Clearly,
\begin{equation}
\labell{eq:varphi-kappa}
\|\varphi^t-\id\|_{C^1}= O(\|\varphi-\id\|_{C^1}) \text{ and }
\|\kappa^t-\id\|_{C^1}= O(\|\varphi-\id\|_{C^1})
\end{equation}
uniformly in $t$.
In particular, $\kappa^t$ is indeed a diffeomorphism,
fixes $p$, and is, moreover, $C^1$-close to the identity when $\varphi$ is
close to $\id$. Furthermore,
\begin{equation}
\labell{eq:A}
\|A_t-I\|=O(\|\varphi-\id\|_{C^1})
\end{equation}
pointwise near $p$ and uniformly in $t$. Hence, $A$ is invertible near
$p$. Since $p$ is an isolated critical point of $F$ by (GF1$^\prime$), $p$ is
also an isolated zero of $X_F$, and thus an isolated zero of
$X_{\tK_t}$. This gives a direct proof of (ii) for $\tK$.

\subsubsection{Property (iii) for $\tK$}
Since
$X_F(p)=0$, the linearization \eqref{eq:XK2} at $p$ yields
$$
d(X_{\tK_t})_p=A_t(p)\circ d(X_F)_p\circ d(\kappa^t)_p.
$$
Here, $d(X_{\tK_t})_p$ and $d(X_F)_p$ are the linear Hamiltonian
vector fields on $T_p\R^{2n}=\R^{2n}$ with quadratic Hamiltonians
$d^2(\tK_t)_p$ and, respectively, $d^2 F_p$.  Furthermore, it is
easy to see from \eqref{eq:A} and \eqref{eq:varphi-kappa} that
$A_t(p)$ and $d(\kappa^t)_p$ are both close to $I$, with error
$O(\|d^2 F_p\|)$, and hence are small. Combining these observations,
we conclude that
$$
\|d^2(\tK_t)_p\|=O(\|d^2 F_p\|)
$$
proving (iii) for $\tK$.

\subsubsection{The upper bound \eqref{eq:KF1} for $\tK$}
\labell{sec:proofKF1}
Turning to the proof of (iv) for $\tK$, observe that for every small
$R>0$ there exists $r>0$ such that $(\bx_t,\by_t)$ and $(\bx_t,y)$ are
in $B_R$ for all $t\in [0,\,1]$ and all $(x,y)$ in $B_r$. Furthermore,
it is clear that
\begin{equation}
\labell{eq:rR}
R=O(\|\varphi-\id\|_{C^1})\cdot r=O_\varphi(r).
\end{equation}

To establish the upper bound \eqref{eq:KF1} of (iv), let us first show
that
\begin{equation}
\labell{eq:KK}
\| X_{\tK_t}(z)-X_{\tK_t}(\kappa^t(z))\|
=\big(O(\|d^2 F_p\|)+O_{\varphi}(r)\big)\cdot\|X_F(\kappa^t(z))\|
\end{equation}
for every $z$ in $B_r$ and all $t\in [0,\,1]$. We have
\begin{equation*}
\begin{split}
\| X_{\tK_t}(z)-X_{\tK_t}(\kappa^t(z))\|
&=\left\|\int_0^1\frac{d}{ds}X_{\tK_{t}}\left(sz+(1-s)\kappa^t(z)\right)\,ds
\right\| \\
&\leq \int_0^1 \|dX_{\tK_{t}}\left(sz+(1-s)\kappa^t(z)\right)\|\,ds
\cdot \|z-\kappa^t(z)\| \\
& \leq \max_{w\in B_R} \|dX_{\tK_{t}}(w)\|\cdot \|z-\kappa^t(z)\| \\
&\leq \max_{w\in B_R} \|dX_{\tK_{t}}(w)\|\cdot \|X_F(\kappa^t(z))\|,
\end{split}
\end{equation*}
where in the last inequality we used the fact that, by \eqref{eq:varphi-F},
$$
\|z-\kappa^t(z)\|=\|y-\by^t\|=\|\p_1 F_t(\kappa^t(z))\|
\leq \|X_F(\kappa^t(z))\|.
$$
Thus, we only need to show that
$$
\max_{B_R} \|dX_{\tK_{t}}\|=O(\|d^2 F_p\|)+O_{\varphi}(r).
$$
By \eqref{eq:XK2}, we have
\begin{equation*}
\begin{split}
\max_{w\in B_R}\|dX_{\tK_t}(w)\|
&\leq
\max_{w\in B_R}
\left\|\left(dA_t(\kappa^t(w)) d\kappa^t(w)\right) X_F(\kappa^t(w)) \right\| \\
&\quad +
\max_{w\in B_R}
\left\|A_t(\kappa^t(w))\left(dX_F(\kappa^t(w))d\kappa^t(w)\right)\right\|
\end{split}
\end{equation*}

Since $\kappa^t(p)=p$ and $X_F(p)=0$, the first summand is obviously
$O_\varphi(r)$. (When $z\in B_r$, both $w=\kappa^t(z)$ and
$\kappa^t(w)$ are, by
\eqref{eq:rR}, in the ball of radius $O_\varphi(r)$.) The second summand is
bounded as
\begin{multline*}
\max_{w\in B_R}
\left\|A_t(\kappa^t(w))\left(dX_F(\kappa^t(w))d\kappa^t(w)\right)\right\| \\
\leq
\max_{w\in B_R}\left\|A_t(\kappa^t(w))\right\| \cdot
\max_{w\in B_R} \left\|dX_F(\kappa^t(w))\right\|\cdot
\max_{w\in B_R} \left\|d\kappa^t(w)\right\|.
\end{multline*}
Here, the first and the last factors are $O(\|\varphi-\id\|_{C^1})$,
and hence bounded from above by a constant
independent of $\varphi$, when $\varphi$ is sufficiently close to
$\id$.  The middle factor is $O(\|d^2 F_p\|)+O_{\varphi}(r)$, for
$\|dX_F(p)\|=\|d^2 F_p\|$ and, as a consequence,
$$
\left\|dX_F(\kappa^t(w))\right\|=O(\|d^2 F_p\|)+O_{\varphi}(\|\kappa^t(w)\|).
$$
Thus, the second summand is $O(\|d^2 F_p\|)+O_{\varphi}(r)$, which completes
the proof of~\eqref{eq:KK}.

Then
\begin{equation*}
\|X_{\tK_t}(\kappa^t(z))-X_F(\kappa^t(z))\|
\leq \|X_{\tK_t}(z)-X_{\tK_t}(\kappa^t(z))\|
+ \|X_{\tK_t}(z)-X_F(\kappa^t(z))\|.
\end{equation*}
By \eqref{eq:XK2}, the second term is bounded as
\begin{equation*}
\begin{split}
\|X_{\tK_t}(z)-X_F(\kappa^t(z))\|
&\leq \| A_t(\kappa^t(z))-I\|\cdot \|X_F(\kappa^t(z))\| \\
 &=\left(O(\|d^2 F_p\|)+O_\varphi(r)\right)\cdot \|X_F(\kappa^t(z))\|
\end{split}
\end{equation*}
and the first term is
$\big(O(\|d^2 F_p\|)+O_{\varphi}(r)\big)\cdot\|X_F(\kappa^t(z))\|$
by \eqref{eq:KK}. This proves the pointwise estimate \eqref{eq:KF1}
at $\kappa^t(z)$ in place of $z$. Since $\kappa^t$ is a diffeomorphism
fixing $p$, the upper bound \eqref{eq:KF1} in its original form (at $z$)
follows.

\subsubsection{The upper bound \eqref{eq:KF2} for $\tK$}
Arguing exactly as in the proof of \eqref{eq:KK}, it is easy to show
that
$$
\| X_{F}(z)-X_{F}(\kappa^t(z))\|
=\big(O(\|d^2 F_p\|)+O_{\varphi}(r)\big)\cdot\|X_F(\kappa^t(z))\|
$$
and, as a consequence,
$$
\big[1-\big(O(\|d^2 F_p\|)+O_{\varphi}(r)\big)\big]\cdot\| X_{F}(z)\|\leq
\|X_F(\kappa^t(z))\|.
$$
Therefore, to establish \eqref{eq:KF2} for $\tK$, it is sufficient to
prove the upper bound
\begin{equation}
\labell{eq:KF22}
\|\dot{X}_{\tK_t}(z)\|=
\big(O(\|d^2 F_p\|)+O_{\varphi}(r)\big)\cdot\|X_F(\kappa^t(z))\|
\end{equation}
for $z\in B_r$.

Differentiating \eqref{eq:XK2} with respect to $t$ and setting
$w=\kappa^t(z)$, we obtain
\begin{equation}
\labell{eq:three-terms}
\begin{split}
\dot{X}_{\tK_t}(z)
&=\dot{A}_t(w)X_F(w)\\
&\quad +
\left(dA_t(w)\dot{\kappa}^t(w)\right)X_F(w)\\
&\quad +
A_t(w)\left(dX_F(w)\dot{\kappa}^t(w)\right).
\end{split}
\end{equation}
To prove \eqref{eq:KF22}, we will estimate all three terms in this identity.
As a straightforward calculation shows,
$$
\dot{A}_t=-(I+\p_{12}F_t)^{-2}
\begin{bmatrix}
\p_{12}F & 0\\
\p_{22}F_t(x,y)\p_{12}F-(I+\p_{12}F_t)\p_{11}F & 0
\end{bmatrix}.
$$
Thus, $\|\dot{A}_t\|=O(\|d^2F_p\|)+O_\varphi(r)$ and
$$
\|\dot{A}_t(w)X_F(w)\|
=\big(O(\|d^2F_p\|)+O_\varphi(r)\big)\cdot \|X_F(w)\|.
$$
Furthermore,
$$
\dot{\kappa}^t(w)=
\left(0, -\big(I-\p_{12}F_t(w)\big)^{-1}\p_1 F(w)\right)
$$
as follows from the definition of $\kappa^t$ and \eqref{eq:varphi-F}.
Using again the inequality $ \|\p_1 F(w)\|\leq \|X_F(w)\|$, we see that
\begin{equation*}
\begin{split}
\|\dot{\kappa}^t(w)\|&\leq O(\|\varphi-\id\|_{C^1})\cdot \|X_F(w)\|\\
&\leq \const \cdot\|X_F(w)\|,
\end{split}
\end{equation*}
where $\const$ is independent of $\varphi$, when $\varphi$ is
sufficiently close to $\id$.  Then, since $\|X_F(w)\|=O_\varphi(r)$,
we have
$$
\left\|\left(dA_t(w)\dot{\kappa}^t(w)\right)X_F(w)\right\|
=O_\varphi(r)\cdot \|X_F(w)\|
$$
and
$$
\left\|A_t(w)\left(dX_F(w)\dot{\kappa}^t(w)
\right)\right\|
=\big(O(\|d^2F_p\|)+O_\varphi(r)\big)\cdot \|X_F(w)\|,
$$
for $\|dX_F(w)\|=O(\|d^2F_p\|)+O_\varphi(r)$.

These estimates combined with \eqref{eq:three-terms} yield
\eqref{eq:KF22}, and hence \eqref{eq:KF2} for $\tK$.

\subsubsection{The Hamiltonian $K$; properties (i), (ii), and (iii)}
Fix a monotone increasing function $\lambda\colon [0,\, 1] \to [0,\,
1]$ such that $\lambda(t)\equiv 0$ when $t$ is near $0$ and
$\lambda(t)\equiv 1$ when $t$ is near $1$. This function is independent of
$\varphi$, and hence $|\lambda'|$ and $|\lambda''|$ are bounded from above by
constants independent of $\varphi$.

The Hamiltonian $K$ is the one generating the flow
\begin{equation}
\labell{eq:K-flow}
\varphi^t_K=\varphi^{t-\lambda(t)}_F\varphi^{\lambda(t)}_{\tK}.
\end{equation}
Explicitly, since $F$ is autonomous,
\begin{equation}
\labell{eq:K}
K_t(z)=(1-\lambda'(t))F(z)+
\lambda'(t)\tK_{\lambda(t)}(\varphi^{\lambda(t)-t}_F(z)).
\end{equation}
It is clear that $K_t\equiv F$  when $t$ is close to $0$ and $1$ and
hence $K$ can be viewed as a Hamiltonian one-periodic in $t$. Also,
$\varphi^1_{K}=\varphi^1_{\tK}=\varphi$, i.e., requirement (i) is
satisfied. As has been pointed out, the second condition, (ii),
follows from (GF1$^\prime$) and (iv) which is proved below.

Passing to the Hessians of the Hamiltonians in \eqref{eq:K} at $p$, we have
$$
d^2(K_t)_p=(1-\lambda'(t))d^2 F_p+
\lambda'(t)d^2(\tK_{\lambda(t)})_p\circ d(\varphi^{\lambda(t)-t}_F)_p.
$$
By (iii) for $\tK$ and \eqref{eq:varphi-kappa},
$\|d^2(K_t)_p\|=O(\|d^2 F_p\|)$, which proves (iii) for $K$.

\subsubsection{The upper bound \eqref{eq:KF1} for $K$}
\labell{eq:KF1K}

The Hamiltonian vector field of $K$ is
\begin{equation}
\labell{eq:XK}
X_{K_t}(z)=(1-\lambda'(t))X_F(z)
+\lambda'(t)d\varphi^{t-\lambda(t)}_F(w)
\big(X_{\tK_{\lambda(t)}}(w)\big),
\end{equation}
where $w=\varphi^{\lambda(t)-t}_F(z)$. Since $F$ is autonomous,
$X_F(z)=d\varphi^{t-\lambda(t)}_F(w) \left(X_F(w)\right)$.
In particular,
\begin{equation}
\labell{eq:zw}
\|X_F(w)\|\leq \const \cdot \|X_F(z)\|,
\end{equation}
when $z$ is close to $p$.  Also note that when $z\in B_r$, the point
$w=\varphi^{\lambda(t)-t}_F(z)$ is in $B_R$, where the radius $R$ satisfies
\eqref{eq:rR}.

With these facts in mind, we have
\begin{equation*}
\begin{split}
\|X_{K_t}(z)-X_F(z)\|&\leq |\lambda'(t)|\left\|
d\varphi^{t-\lambda(t)}_F(w) \left(X_{\tK_{\lambda(t)}}(w)\right)
-X_F(z)\right\|\\
&\leq |\lambda'(t)|\left\|
d\varphi^{t-\lambda(t)}_F(w) \left(X_{\tK_{\lambda(t)}}(w)\right)
-d\varphi^{t-\lambda(t)}_F(w) \left(X_F(w)\right)\right\|\\
&\leq |\lambda'(t)|\big\|d\varphi^{t-\lambda(t)}_F(w)\big\|\cdot
\left\|X_{\tK_{\lambda(t)}}(w)-X_F(w)\right\|\\
&= \big( O(\|d^2 F_p\|)+O_\varphi(r)\big)\|X_F(w)\|\\
&= \big( O(\|d^2 F_p\|)+O_\varphi(r)\big)\|X_F(z)\|,
\end{split}
\end{equation*}
when $z\in B_r$. Here, the next to the last estimate follows from
\eqref{eq:KF1} for $\tK$. This proves \eqref{eq:KF1} for $K$.

\subsubsection{The upper bound \eqref{eq:KF2} for $K$}
Differentiating \eqref{eq:XK} with respect to $t$, we obtain
\begin{equation*}
\begin{split}
\dot{X}_{K_t}(z)&=\lambda''(t)\left(
d\varphi^{t-\lambda(t)}_F(w)\big(X_{\tK_{\lambda(t)}}(w)\big)-X_F(z)\right)\\
&\quad+ \lambda'(t)^2 d\varphi^{t-\lambda(t)}_F(w)
\big(\dot{X}_{\tK_{\lambda(t)}}(w)\big)\\
&\quad+ \lambda'(t)(\lambda'(t)-1)
d\varphi^{t-\lambda(t)}_F(w)[X_F,X_{\tK_t}](w).
\end{split}
\end{equation*}
Arguing as in Section \ref{eq:KF1K}, we see that the norm of the
first term in this sum is $\big( O(\|d^2
F_p\|)+O_\varphi(r)\big)\|X_F(z)\|$. Similarly, the same holds for
the second term by \eqref{eq:KF2} for $\tK$. (In both cases we use
\eqref{eq:zw} to relate $X_F(z)$ and $X_F(w)$ and also the fact that
$\|d\varphi^{t-\lambda(t)}_F(w)\|\leq \const$ when $z\in B_r$.)

To estimate the third term,  it is sufficient to show that
\begin{equation}
\labell{eq:bracket}
\|[X_F,X_{\tK_t}](w)\|=\big( O(\|d^2 F_p\|)+O_\varphi(r)\big)\|X_F(w)\|,
\end{equation}
for then, by \eqref{eq:zw}, this term is $\big( O(\|d^2
F_p\|)+O_\varphi(r)\big)\|X_F(z)\|$ when $z\in B_r$.
To prove \eqref{eq:bracket}, observe that
\begin{equation*}
\begin{split}
\big\|[X_F,X_{\tK_t}](w)\big\|&=\big\|[X_F,X_{\tK_t}-X_F](w)\big\|\\
&\leq \alpha(w)\|X_F(w)\|+\beta(w)\|X_{\tK_t}(w)-X_F(w)\|,
\end{split}
\end{equation*}
where the functions $\alpha(w)\geq 0$ and $\beta(w)\geq 0$ are bounded
from above by the partial derivatives of $X_{\tK_t}-X_F$ and,
respectively, $X_F$ at $w$. Hence, both of these functions are
$O(\|d^2 F_p\|)+O_\varphi(r)$ and \eqref{eq:bracket} follows from
\eqref{eq:KF1} for $\tK$.

This completes the proof of \eqref{eq:KF2} for $K$ and the proof of the lemma.

\subsection{Proof of Lemma \ref{lemma:special}}
Let $\Phi=d\varphi_p$ and let $V=L\oplus L'$ be the decomposition of
$V=T_pW$ from Lemma \ref{lemma:la2}. Pick a linear canonical
coordinate system $(x,y)$ on $T_pW$, which is compatible with the
decomposition, i.e., such that the $x$-coordinates span $L$ and the
$y$-coordinates span $L'$. By Lemma \ref{lemma:la2}, we can do this so
that $\|\Phi-I\|$ is small in this coordinate system, and thus
$\varphi$ is given by the generating function $F$.  Denote by $Q$ the
Hessian of $F$ at $p$ and let $X_Q$ be the linear Hamiltonian vector
field of $Q$ on $V$.

Linearizing \eqref{eq:varphi-F} at $p$, we see that $\Phi$ and
$Q$ are related via the equation
\begin{equation}
\labell{eq:Phi-Q}
\Phi-I=X_Q P(\Phi).
\end{equation}
Here $P(\Phi)\colon V\to V$ is obtained from $\Phi$ by replacing its
$y$-component by the identity map, i.e.,
$P(\Phi)(x,y)=(\bx,y)$ in the decomposition $V=L\oplus L'$, where
$\Phi(x,y)=(\bx,\by)$. Note that \eqref{eq:Phi-Q} uniquely determines
$X_Q$.

Furthermore, let $\tPhi_t$ be the linearization of $\varphi_{\tK}^t$ at
$p$. This family of linear symplectic transformations satisfies the
equation
\begin{equation}
\labell{eq:tPhi-Q}
\tPhi_t=I+tX_Q P(\tPhi_t),
\end{equation}
which again uniquely determines $\tPhi_t$.

It is clear from \eqref{eq:Phi-Q} and \eqref{eq:tPhi-Q} that $X_Q$ and
$\tPhi_t$ depend only on the decomposition $V=L\oplus L'$. Hence, any other
coordinate system compatible with this decomposition will give rise to the
same quadratic form $Q$ and the same maps $\tPhi_t$.

Due to \eqref{eq:K-flow}, the linearization $d(\varphi_K^t)_p$
is equal to
$$
\Phi_t=\exp\big( (t-\lambda(t)) X_Q\big) \tPhi_{\lambda(t)}.
$$
Hence, $\Phi_t$ also depends only on the decomposition, but not on
the coordinate system as long as the latter is compatible with the
decomposition. In other words, all such coordinate systems result in
the flows $\varphi^t_K$ with linearization $\Phi_t$.

Lemma \ref{lemma:special} follows now from Lemma \ref{lemma:la2},
which guarantees that there exist symplectic bases $\Xi$ (or,
equivalently, linear canonical coordinate systems) compatible with
$V=L\oplus L'$ and making $\|\Phi-I\|_{\Xi}$ arbitrarily small.

\begin{Remark}
Recall that $Q\leq 0$ due to (K1) and that all eigenvalues of $\Phi$
are equal to one. Combining these facts with the normal forms of
quadratic Hamiltonians (see \cite[Appendix 7]{Ar} and \cite{Wi}), it
is not hard to show that $Q=-(y_1^2+\ldots+y_k^2)$ in some symplectic
basis compatible with the decomposition $T_pW=L\oplus L'$.
Then it is straightforward to write down an explicit expression for
$\tPhi_t$ and $\Phi_t$. This, however, does not lead to any
simplification in the line of reasoning used here, for the
required result readily follows from \eqref{eq:Phi-Q} and
\eqref{eq:tPhi-Q}.
\end{Remark}

\section{Proof of Proposition \ref{prop2}}
\labell{sec:prop2-pf}

\subsection{Outline of the proof}
\labell{sec:prop2-outline}
First note that it is sufficient to prove the proposition for
the Hamiltonian $K^1$ in place of $H$ and the constant orbit $p$
of $K^1$ in place of $\gamma$.

Indeed, $p$ is a symplectically degenerate maximum of $K^1$ as is
pointed out in Example \ref{ex:sympl-deg}.  The Hamiltonians $K^1$ and
$H$ have the same time-$T$ flow and there is a natural one-to-one
correspondence between (contractible) $T$-periodic orbits of the
Hamiltonians, for $\varphi^t_{H}=\eta^t\circ\varphi^t_{K^1}$ with $t\in
S^1_T$. Due to our normalization of $K^1$, the corresponding
$T$-periodic orbits of $K^1$ and $H$ have equal actions and, in
particular, $(K^1)^{(T)}$ has the same action spectrum and action
filtration in the Floer complex as $H^{(T)}$; see
Section~\ref{sec:loops}. As a consequence,
$$
\HF^{(Tc+\delta_T,\,Tc+\eps)}_{n+1}\big(H^{(T)}\big)=
\HF^{(Tc+\delta_T,\,Tc+\eps)}_{n+1}\big((K^1)^{(T)}\big).
$$
Thus, the proposition holds for $H$ if (and only if) it holds for $K^1$.
Furthermore,
when $H$ is replaced by $K^1$, the loops $\eta_i^t$ get replaced by the loops
$\eta_i^t\circ(\eta_1^t)^{-1}$ which have the identity linearization at $p$
by (K3).

To summarize, keeping the notation $H$ for the Hamiltonian $K^1$, we
may assume throughout the proof that the
Hamiltonian $H$ is such that
\begin{itemize}

\item the point $p$ is a strict local maximum of $H_t$ for all
$t\in S^1$, and
\item $d(\eta_i^t)_p=I$ for all $t\in S^1$.
\end{itemize}

With these observations in mind, we establish the proposition by
using the squeezing method of \cite{BPS,GG}.  Namely, closely
following \cite{GG}, we construct functions $H_\pm$ such that
$H_-\leq H\leq H_+$ (see Fig.\ 1) and such that the map $\Psi_{H_+,H_-}$
in the filtered Floer
homology for the interval $(Tc+\delta_T,\,Tc+\eps)$ induced by a
monotone homotopy from $H_+$ to $H_-$ is non-zero. This map factors
as
$$
\HF_{n+1}^{(Tc+\delta_T,Tc+\eps)}\big(H_+^{(T)}\big)\to
\HF_{n+1}^{(Tc+\delta_T,\,Tc+\eps)}\big(H^{(T)}\big)\to
\HF_{n+1}^{(Tc+\delta_T,\,Tc+\eps)}\big(H_-^{(T)}\big),
$$
and, therefore, $\HF_{n+1}^{(Tc+\delta_T,\,Tc+\eps)}\big(H^{(T)}\big)\neq 0$
as required.

The Hamiltonian $H_+$ depends only on $H$ and $\eps$.  Outside a ball $B_R$ of
radius $R>0$, centered at $p$, the function $H_+$ is constant and
equal to $\max H_+$.  (Here the distance is taken with respect to some
fixed metric compatible with $\omega$.)  Within $B_R$, the Hamiltonian
$H_+$ is a function of the distance to $p$, equal to $c=H(p)$ when the
distance is small, dropping to some constant $a<c$, and then
increasing to $\max H_+$ near the boundary of $B_R$.

\begin{figure}[ht!] 
\labellist 
\small\hair 2pt 
\pinlabel $H_+$ at 175 117
\pinlabel $H$ at 174 82 
\pinlabel $H_-$ at 175 40
\pinlabel $c$ [r] at 0 85 
\pinlabel $a$ [r] at 0 75 
\pinlabel $p$ [bl] at 0 -11
\pinlabel $W$ [br] at 198 -11 
\endlabellist 
\centering 
\includegraphics[scale=0.8]{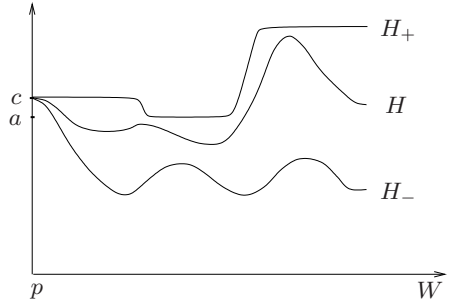} 
\caption{The functions $H$ and $H_\pm$}
\end{figure}

The period $T$ is required to be large enough, i.e., $T\geq T_0$, where $T_0$
is determined by $H_+$. The function $H_-$ and
the constant $\delta_T>0$ depend on $T\geq T_0$. The condition that $p$
is a symplectically degenerate maximum of $H$ is used in the
construction of $H_-$ and also in proving that $\Psi_{H_+,H_-}\neq 0$.

The function $H_-$ is constructed as follows.  Pick $i$ so that
$T\cdot\| d^2 (K_t^i)_p\|_{\Xi^i}$ is small.  (Here, as above, $K^t_i$
is normalized by $K_t^i(p)\equiv c$.)  There exists a bump function
$F\leq K_i$, supported near $p$, with non-degenerate maximum at $p$
and $F(p)=c$ and such that $T\cdot\| d^2 F_p\|_{\Xi^i}$ is also
small. Then $\Psi_{H_+,F}\neq 0$. Setting $H_-$ to be the Hamiltonian
generating the flow $\eta_i^t\circ\varphi_F^t$, normalized by
$H_-(p)\equiv c$, we note that $H_-\leq
H$. Hence, $H_-\leq H\leq H_+$.  The Hamiltonian $H_-$ has the same
filtered Floer homology as $F$, and we show that $\Psi_{H_+,H_-}\neq
0$ using (K3).

\subsection{Bump functions}
\labell{sec:bump}
In this section we recall a few standards facts needed in the proof,
concerning the filtered Floer complex of a bump function.

\subsubsection{Bump functions on $\R^{2n}$}
\labell{sec:bump-r2n}

In the standard canonical coordinates $(x,y)$ on $\R^{2n}$, set
$\rho=(x^2+y^2)/2$. All orbits of $\varphi^t_\rho$ are
closed and have period $2\pi$.
Fix a ball $B_{r}\subset \R^{2n}$ of radius $r>0$, centered at
the origin $p$.

Consider a rotationally symmetric function $F$ on $\R^{2n}$ supported
in $B_{r}$. The function $F$ depends only on the distance to $p$
and it will be convenient in our analysis to also view $F$ as a
function of $\rho$.  Assume, in addition, that $F$ has the following
properties (see Fig.\ 2):

\begin{itemize}

\item $F$ is decreasing as a function of $\rho$;

\item $|F'|<\pi$ and $|F'|$ is increasing
on some closed ball $\bar{B}_{r_-}\subset B_{r}$; and

\item on the shell $B_{r}\ssminus B_{r_-}$,

\begin{itemize}
\item[$\circ$] $F$ is concave, i.e.,
$F''\leq 0$, on
$[r^2_-/2,\, (r')^2/2]$, where $r_-<r'<r$,

\item[$\circ$] $F$ is convex ($F''\geq
0$) on $[(r'')^2/2,\, r^2/2]$, where $r'<r''<r$, 

\item[$\circ$] $F$ has constant slope
($F'= \const$) on the interval $[(r')^2/2,\, (r'')^2/2]$, where
$\const/\pi$ is irrational.
\end{itemize}
\end{itemize}
We will refer to $F$ as a \emph{standard bump function} on
$\R^{2n}$ and  we will call $C:=F(p)$ and $r_-$ and
$r$ and other constants from the construction of $F$ the
parameters of $F$.

\begin{figure}[ht!] 
\labellist 
 \small\hair 2pt 
 \pinlabel $F$ at -3 79 
 \pinlabel $C$  at -3 59
\pinlabel $p$   at 2 -4
\pinlabel $x$   at 40 -4
\pinlabel $y$   at 56 -4
\pinlabel $\rho$ at 102 -4
\endlabellist 
\centering 
\includegraphics[scale=1.2]{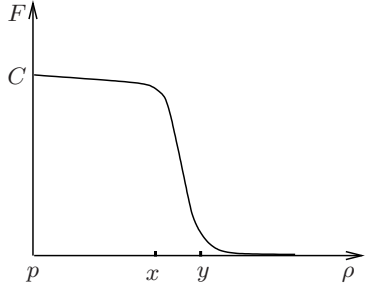} 
\caption{The bump function $F$}
\end{figure}

The trivial one-periodic orbits of $F$ (i.e., its critical points) are
either contained in $\bar{B}_{r'}$ or in the complement to $B_{r''}$.
The orbits from the first group form a closed ball (possibly of zero
radius) centered at $p$ and have action $C$; the orbits from the
second group are exactly the points where $F=0$.

Non-trivial one-periodic orbits fill in spheres of radii $r^\pm_i$ with
$$
r_-< r_1^-< r_2^- < \ldots  <r'\text{ and }
r''<\ldots< r_2^+ < r_1^+ < r.
$$

Let $\tF$ be the standard $C^2$-small (periodic in
time) perturbation of $F$ as the ones considered
in, e.g., \cite{FHW,GG}, and still supported in $B_{r}$.
For such a perturbation each sphere filled in by one-periodic orbits of
$F$ breaks down into $2n$ non-degenerate orbits. Within $\bar{B}_{r_-}$,
we may assume that $\tF$ is still autonomous and rotationally
symmetric and $0<|\tF'|<\pi$. (Hence the only one-periodic orbit of $\tF$
in $\bar{B}_{r_-}$ is the trivial orbit $p$ of Conley-Zehnder index $n$.)

As is well known, the filtered Floer complex of $\tF$ and the filtered
Floer homology of $\tF$ (and $F$) are still defined, say, for any
positive interval of actions $0\leq a<b$ even though $\R^{2n}$ is not
compact; see, e.g., \cite{FH,FHW,FS,CGK,Gi,GG}. Here we adopt the
conventions of \cite{CGK,Gi,GG}.

Pick $\eps>\delta>0$ so that
\begin{equation}
\labell{eq:eps-del}
\eps>\pi r^2,\text{ and } C-\delta>2\pi r^2 \text{ and } \delta<\pi r_-^2.
\end{equation}
Assuming that $C>\pi r^2$, we are interested in the periodic orbits
of $\tF$ with indices $n-1$ or $n$ or $n+1$ and action in the range
$(C-\delta,\, C+\eps)$.

It is not hard to see that $\tF$ has no one-periodic orbits of index
$n-1$. Furthermore, it has exactly two one-periodic orbits
of index $n$. One of these is the constant orbit $p$. The second
orbit $y$ arises from the sphere of
periodic orbits farthest from the origin. This sphere has radius
$r_1^+$ and $A_{\tF}(y)=\pi (r_1^+)^2$, up to an error of order
$\|F-\tF\|_{C^1}$. By \eqref{eq:eps-del}, $y$ is outside the range of
action.  Finally, $\tF$ has two one-periodic orbits of index $n+1$,
but only one of them, $x$, has action in $(C-\delta,\, C+\eps)$.  (The
action of the second orbit is approximately equal to $\pi (r_1^+)^2$.)
The orbit $x$ arises from the sphere of periodic orbits closest to the
origin.  This is the sphere of radius $r_1^-$ and $A_{\tF}(x)=C+\pi
(r_1^-)^2$, up to an error of order $\|F-\tF\|_{C^1}$.  (See \cite{CGK,GG} or
Section \ref{sec:K+}, where we analyze in detail periodic orbits of a
function similar to~$F$.)

It is well known that $\p x=y+p$ in the Floer complex of $\tF$ for
the interval $(0,\,\infty)$; see, e.g., \cite{Gi,GG}.
Summarizing these observations, we conclude that

\begin{itemize}

\item $\tF$ has no one-periodic orbits of index $n-1$;

\item $p$ is the only one-periodic orbit of $\tF$ with index $n$ and
action in $(C-\delta,\, C+\eps)$;

\item $x$ is the only one-periodic orbit of $\tF$ with index $n+1$ and
action in $(C-\delta,\, C+\eps)$;

\item the connecting map from the long exact sequence
$$
\Z_2\cong
\HF_{n+1}^{(C+\delta,\, C+\eps)}(F)\to \HF_{n}^{(C-\delta,\, C+\delta)}(F)
\cong\Z_2
$$
is an isomorphism sending $[x]$ to $[p]$, and hence
$\HF_{n}^{(C-\delta,\, C+\eps)}(F)=0$.

\end{itemize}

\subsubsection{Bump functions on a closed manifold}
\labell{sec:bump-W}

Let $U$ be a small neighborhood of $p\in W$. Fixing a canonical
coordinate system on $U$, denote the open ball of radius $r>0$ in $U$,
centered at $p$, by $B_r$ and let $S_r$ be the boundary of $B_r$.

We define a bump function $F$ on $W$ exactly as for $\R^{2n}$ by using
the coordinate system in $U$. Furthermore, since $W$ is compact, now we
need not assume that $F$ is supported in $B_{r}$. Instead we just
require $F$ to be constant outside $B_{r}$. In other words, we allow
$F$ to be shifted up and down.

The description of periodic orbits of $\tF$ and the Floer homology of
$F$ given in Section \ref{sec:bump-r2n} extends word-for-word to this
case, provided that $B_{r}$ is sufficiently small (e.g., displaceable)
and the variation $C-\min F$ is sufficiently large. The requirement
$C>\pi r^2$ is replaced by that $C-\min F> h(B_{r})$, where $h(B_{r})$
depends only on $B_{r}$ and goes to zero as $r\to 0$; see, e.g.,
\cite{Gi}. (Hypothetically, $h(B_{r})$ is equal to the displacement
energy of $B_{r}$, although the estimate we have been able to prove is
somewhat weaker.) The requirement \eqref{eq:eps-del} carries over to
this case unchanged when $F$ is supported in $B_{r}$, and is, in
general, replaced by
\begin{equation}
\labell{eq:eps-del2}
\eps>\pi r^2, \text{ and } C-\delta>\min F+2\pi r^2
\text{ and } \delta<\pi r_-^2.
\end{equation}

\subsubsection{Connecting trajectories from $x$ to $p$}
\labell{sec:traj} Let us show that by making $r$ sufficiently small,
we can ensure that the Floer gradient trajectories of $F$ from $x$ to
$p$ are close to $p$. (When $F$ is a bump function on $\R^{2n}$, all such
trajectories are contained in $\bar{B}_r$ by the maximum principle.)

To this end, pick a ball $B_R\supset \bar{B}_{r}$ contained in $U$ and
fix once and forever a compatible with $\omega$ almost complex
structure $J_0$ on $W$ coinciding with the standard complex structure
on a neighborhood of $\bar{B}_{R}$. Consider holomorphic curves $v$ in
$\bar{B}_{R}\ssminus B_{r}$ with boundary in $S_{r}\cup S_R$ and such
that the part of the boundary of $v$ lying in $S_{r}$ is
non-empty. (Then the part of the boundary of $v$ in $S_R$ is also
non-empty due to the maximum principle.) Denote by $A(r,R)>0$ the
infimum of the areas of such curves $v$. It is easy to see that
$A(r,R)$ remains separated from zero as $r\to 0$. (Otherwise we would
have $A(r,R)=0$ for some fixed $r>0$ as is clear from considering the
intersections with $\bar{B}_{R}\ssminus B_{r}$ of holomorphic curves
whose areas approach zero.)  In other words, $\liminf_{r\to
0+}A(r,R)>0$. Replacing $R$ by $R/2$, we see that there exists
$r_0(R,J_0)>0$ such that
\begin{equation}
\labell{eq:r0}
\pi r^2 < A(r,R/2) \text{ for all positive $r<r_0(R,J_0)<R/2$.}
\end{equation}

\begin{Lemma}
\labell{lemma:traj}

Let $F$ be an arbitrary bump function $F$  such that
\eqref{eq:r0} holds and $C-\min F>h(B_{r})$. Assume that
$\eps>0$ and $\delta>0$ satisfy \eqref{eq:eps-del2}.  Then for a
perturbation $\tF$ of $F$ as above and any regular perturbation of $J$
of $J_0$ all Floer anti-gradient trajectories from $x$ to $p$ are
contained in ${B}_{R}$.
\end{Lemma}

\begin{proof}
Assume the contrary. Then for some $\tF$ close to $F$ and for some sequence
of regular perturbations $J_l\to J_0$, there exists a sequence of
connecting trajectories $u_l$ from $x$ to $p$, leaving a neighborhood
of $\bar{B}_{R/2}$. Observe that the part of $u_l$ contained in
$\bar{B}_{R/2}\ssminus B_{r}$ is a $J_l$-holomorphic curve.  By the
compactness theorem, in the limit we have a $J_0$-holomorphic curve
$v$ in $\bar{B}_{R/2}\ssminus B_{r}$ with non-empty boundary in $S_{r}$.
By the definition of $A(r,R/2)$, the area of $v$ is greater than
$A(r,R/2)$. Therefore, the same is true for the part of $u_l$ contained in
$\bar{B}_{R/2}\ssminus B_{r}$ when $J_l$ is close to $J_0$.  Thus,
$E(u_l)> A(r,R/2)>\pi r^2$ by \eqref{eq:r0}. This is impossible, for
$$
E(u_l)=A_{\tF}(x)-A_{\tF}(p)=\pi (r_1^-)^2+
O\big(\|F-\tF\|_{C^1}\big) <\pi r^2.
$$
\end{proof}

\subsection{The function $H_+$}
\labell{sec:K+}
Without loss of generality, we may assume that $H\geq 0$.
Furthermore, throughout this section we will keep the notation and
convention of Section \ref{sec:bump}. In particular, we fix a system
of canonical coordinates on a neighborhood $U$ of $p$ and let, as in
Section \ref{sec:bump}, the function $\rho$ on $U$ be one half of the
square of the distance to $p$ with respect to this coordinate system.

\subsubsection{The description of $H_+$}
Pick four balls centered at $p$ in $U$:
$$
B_{r_-}\subset B_{r}\subset B_R \subset B_{R_+} \Subset U.
$$

Let $H_+$ be a function of $\rho$, also treated as a function on $U$,
with the following properties (see Fig.\ 1):

\begin{itemize}

\item $H_+\geq H$;

\item $H_+|_{B_{r_-}}\equiv c=H(p)$;

\item on the shell $B_{r}\ssminus B_{r_-}$ the function $H_+$ is monotone
decreasing, as a function of $\rho$, and

\begin{itemize}
\item[$\circ$] $H_+$ is concave ($H_+''\leq 0$) on
$[r_-^2/2,\, (r')^2/2]$, where $r_-<r'<r$, 

\item[$\circ$] $H_+$ is convex ($H_+''\geq
0$) on $[(r'')^2/2,\, r^2/2]$, where $r'<r''<r$, 

\item[$\circ$] $H_+$ has constant slope
($H_+'= \const$) on the interval $[(r')^2/2,\, (r'')^/2]$, where
$\const/\pi$ is irrational;
\end{itemize}

\item $H_+\equiv a$ on the shell $B_R\ssminus B_{r}$, where the constant
$a$ is to be specified latter;

\item $H_+$ is monotone increasing on the shell $B_{R_+}\ssminus B_R$,
and

\begin{itemize}
\item[$\circ$] $H_+$ is convex ($H_+''\geq
0$) on $[R^2/2,\, (R')^2/2]$, where $R<R'<R_+$, 

\item[$\circ$] $H_+$ is concave ($H_+''\leq 0$) on
$[(R'')^2/2,\, R_+^2/2]$, where $R'<R''<R_+$, 

\item[$\circ$] $H_+$, as a function of $\rho$, has constant slope
($H_+'= \const$) on the interval $[(R')^2/2,\, (R'')^2/2]$, where
$\const/\pi$ is irrational;
\end{itemize}

\item $H_+\equiv \max H_+$ on $U\ssminus B_{R_+}$, with the
constant $\max H_+>c$ to be specified.

\end{itemize}

Furthermore, we extend $H_+$ to $W$ by setting it to be constant and equal to
$\max H_+$ on the
complement of $U$. The constant $\max H_+$ is chosen so that $H\leq
H_+$ on $W$ and $\max H_+>c$.

Note that within $B_R$, the function $H_+$ is a standard bump function
of Section~\ref{sec:bump}. This bump function has variation $c-a$ which,
due to the requirement $H_+>H$, may be very small.

\subsubsection{The parameters of $H_+$}
Let us now specify the parameters of
$H_+$. The main, but not the only, requirement on $H_+$ is that $H_+\geq H$.

The neighborhood $U$ is chosen so that $c=H(p)$ is a strict global
maximum of $H$ on $U$ and $U$ is displaceable in $W$ by a
Hamiltonian diffeomorphism. The values $R_+>R>0$ are chosen
arbitrarily, with the only restriction that $B_{R_+}\Subset U$.

To pick $r$, we fix $\eps>0$ and also fix a compatible with
$\omega$ almost complex structure $J_0$ on $W$ coinciding with
the standard complex structure on a neighborhood of the closed ball
$\bar{B}_{R}$. The radius $r$ is chosen so that
\begin{equation}
\labell{eq:r02}
0<r<r_0(R,J_0) \text{ and } \pi r^2 <\eps,
\end{equation}
where the upper bound $r_0(R,J_0)$ is as in Section \ref{sec:traj}.
The radius $r_->0$ is chosen arbitrarily with the only restriction
that $0<r_-<r$.

The constants $a$ and $\max H_+$ are picked so that
$a<c$ and $H_+\geq H$ on $[r_-^2/2,\, R^2/2]$ and
on $W\ssminus B_{R_+}$ (This may require $a$ to be very
close to $c$.) Likewise, on the intervals $[r_-^2/2,\,r^2/2]$ and
$[R^2/2,\,R_+^2/2]$, the behavior of $H_+$ is specified to guarantee that
$H_+\geq H$. Finally, we will also impose the condition
that
\begin{equation}
\labell{eq:c-eps}
\max H_+>c+\eps.
\end{equation}

At this stage we fix $H_+$ satisfying the above requirements.

\subsubsection{Periodic orbits of $H_+^{(T)}$}
\labell{sec:orbits} In this section, we analyze the relevant
$T$-periodic orbits of $H_+$ when $T$ is sufficiently large.  Since $H_+$ is
autonomous, its $T$-periodic orbits can simply be treated as
one-periodic orbits of $T\cdot H_+$. Furthermore, it is clear that all
$T$-periodic orbits of $H_+$ outside $U$ are trivial. Those in $U$ are
either trivial or fill in spheres of certain radii. Replacing $H_+$ by
its standard time-dependent $C^2$-small perturbation $\tH_+$ as in
\cite{CGK,FHW,GG} and Section \ref{sec:bump} results in each of these
spheres splitting into $2n$ non-degenerate orbits. Here, as in Section
\ref{sec:bump}, we are primarily interested in the orbits with index
$n-1$ or $n$ or $n+1$ and action in the interval
$(Tc-\delta,\,Tc+\eps)$ for some small $\delta>0$. We will show that
these orbits are essentially the same as for a bump function with
(large) variation $T\cdot (c-a)>\pi r^2$.

The perturbation $\tH_+$ is similar to $\tF$ from Section
\ref{sec:bump}. We emphasize that $H_+$ is perturbed not only within
the shells $\bar{B}_{R_+}\ssminus B_{R}$ and $\bar{B}_{r}\ssminus
B_{r_-}$ where non-trivial periodic orbits are, but also within the
ball $B_{r_-}$ where $H_+\equiv c$. On this ball, $\tH_+$ is a
monotone decreasing function of $\rho$ with a non-degenerate maximum
at $p$ equal to $c$. This function is $C^2$-close to the constant
function $H_+$ so that (for a fixed $T$), the function $T\cdot\tH_+$
is $C^2$-close to $Tc$ on $B_{r_-}$. In particular, the eigenvalues of
the Hessian $d^2(T\cdot\tH_+)_p$ are close to zero and the
Conley--Zehnder index of the constant $T$-periodic orbit $p$ of
$\tH_+$ is $n$.  In what follows, we will always assume that $\tH_+$
is as close to $H_+$ as necessary.  In the shell $B_R\ssminus B_{r}$
and in the complement to $B_{R_+}$ we keep $\tH_+$ constant and equal
to $H_+$.

With $H_+$ and $\eps>0$ fixed, assume throughout this subsection that
$T$ is sufficiently large and $\delta>0$ is small or, more
specifically, that
\begin{equation}
\labell{eq:eps-del0}
T\cdot (c-a)>2\pi r^2+\delta, \text{ where } \delta<\pi r_-^2\text{ and }
\delta< c-a.
\end{equation}

The trivial $T$-periodic orbits of $H_+$ are the points of
$\bar{B}_{r_-}$ (with action $Tc$), the points of $\bar{B}_{R}\ssminus
B_{r}$ (with action $Ta$), and the points of $W\ssminus B_R$ (with
action $T\cdot\max H_+$).  Here, only the points of $\bar{B}_{r_-}$
have action within the range in question. Indeed, $T\cdot\max
H_+>Tc+\eps$ by \eqref{eq:c-eps} and $Ta<Tc-\delta$, for
$0<\delta<c-a$ by \eqref{eq:eps-del0}. Thus, $p$ is the only trivial
$T$-periodic orbit of $\tH_+$ with action in $(Tc-\delta,\,Tc+\eps)$;
it has index $n$.

We divide non-trivial $T$-periodic orbits of $H_+$ and $\tH_+$ into
four groups.

The first group is formed by the $T$-periodic orbits in the shell
$\bar{B}_{r'}\ssminus B_{r_-}$. These orbits fill in a finite
number of spheres $S_{r_l^-}$ of radii $r_l^-$ with
$$
r_-< r_1^-< r_2^- < \ldots  <r'.
$$
The orbits on $S_{r_l^-}$ have action $T c +\pi (r_l^-)^2\cdot l$.
Once the Hamiltonian $H_+$ is replaced by $\tH_+$, a sphere
$S_{r_l^-}$ breaks down into $2n$ non-degenerate orbits. The
Conley--Zehnder indices of these orbits  are $(2l-1)n+1,\ldots, (2l+1)n$
as is proved in \cite{CGK,GG}.  Only one of the
orbits in this group has index from $n-1$ to $n+1$ and action in the
range $(Tc-\delta,\,Tc+\eps)$. This is a periodic orbit,
denoted by $x$, of index $n+1$ and action $T c +\pi (r_1^-)^2$ (up to
an error of order $\|H_+-\tH_+\|_{C^1}$), arising from the sphere
$S_{r_1^-}$. The Conley--Zehnder indices of the remaining orbits are
greater than $n+1$ although some of these orbits may have action
within the range $(Tc-\delta,\,Tc+\eps)$.

The second group is comprised of the $T$-periodic orbits in the shell
$\bar{B}_{r}\ssminus B_{r''}$. These orbits fill in the
spheres $S_{r_l^+}$ of radii $r_l^+$ with
$$
r''<\ldots< r_2^+ < r_1^+ <r .
$$
The orbits on $S_{r_l^+}$ have action $T a +\pi (r_l^+)^2\cdot l$.
Again, once the Hamiltonian $H_+$ is replaced by $\tH_+$, each sphere
$S_{r_l^+}$ breaks down into $2n$ non-degenerate orbits. The
Conley--Zehnder indices of these orbits are $(2l-1)n,\ldots,
(2l+1)n-1$; see \cite{CGK,GG}. Only the orbits arising from
$S_{r_1^+}$ and $S_{r_2^+}$ can have index $n-1$ or $n$ or
$n+1$. (Other spheres give rise to orbits of index greater than
$5n-1>n+1$.) However, the orbits coming from the spheres $S_{r_1^+}$
and $S_{r_2^+}$ have action not exceeding $T a +2\pi r^2$ if $\tH_+$
is close to $H_+$. By \eqref{eq:eps-del0}, these orbits are
outside the action range $(Tc-\delta,\,Tc+\eps)$.

The $T$-periodic orbits in the shell $\bar{B}_{R'}\ssminus B_{R}$ are
in the third (possibly empty) group. These orbits fill in the spheres
$S_{R_l^-}$ of radii $R_l^-$ with
$$
R< R_1^-< R_2^- < \ldots  <R',
$$
and the orbits on $S_{R_l^-}$ have action $T a -\pi (R_l^-)^2\cdot
l$. Hence, all of these orbits are outside of the range of action
$(Tc-\delta,\,Tc+\eps)$.

The fourth group, which may also be empty, is
formed by the $T$-periodic orbits in the shell $\bar{B}_{R_+}\ssminus
B_{R''}$. These orbits fill in a finite collection of spheres
$S_{R_l^+}$ of radii $R_l^+$ such that
$$
R''<\ldots< R_2^+ < R_1^+ < R_+,
$$
and the orbits on $S_{R_l^+}$ have action $T\cdot \max H_+-\pi
(R_l^+)^2\cdot l$, which can be in the
interval $(Tc-\delta,\,Tc+\eps)$. However, calculating the
Conley--Zehnder indices of the resulting orbits of $\tH_+$ as in
\cite{CGK,GG}, it is easy to see that the sphere $S_{R_l^+}$ brakes down
into non-degenerate orbits of $\tH_+$ with indices $-(2l+1)n+1,\ldots,
-(2l-1)n$. In particular, all resulting orbits have indices not
exceeding $-n$, and none of the orbits has index $n-1$ or $n$ or
$n+1$.

To summarize, the perturbation $\tH_+$ has only one $T$-periodic orbit
of index $n$ with action in $(Tc-\delta,\,Tc+\eps)$ -- this is the
trivial orbit $p$ -- and only one orbit, namely $x$, of index $n+1$
with action within this range. The action of $p$ is $Tc$ and the
action of $x$ is $T c +\pi (r_1^-)^2$ up to an error of order
$\|H_+-\tH_+\|_{C^1}$. There are no orbits with index $n-1$ and action in
the range $(Tc-\delta,\,Tc+\eps)$.

\subsubsection{The Floer homology of $H_+^{(T)}$}
\labell{sec:FH-K+}
As in the previous section, assume that $T$ is sufficiently large
and $\delta>0$ is small (independently of $T$).
Explicitly, now we require in addition to \eqref{eq:eps-del0} that
\begin{equation}
\labell{eq:T-delta}
T(c-a)> h(B_{R}),
\end{equation}
where $h(B_{R})$ is defined in Section \ref{sec:bump-W}. In this
section we prove

\begin{Lemma}
\labell{lemma:FH1}
Under the above assumptions on the function $H_+$, the period $T$, and
$\eps$ and $\delta$, we have
$$
\HF_n^{(Tc-\delta,\,Tc+\delta)}\big(H_+^{(T)}\big)\cong\Z_2 \text{ and }
\HF_{n+1}^{(Tc+\delta,\,Tc+\eps)}\big(H_+^{(T)}\big)\cong\Z_2
$$
with generators $[p]$ and, respectively, $[x]$.
Moreover, the connecting map 
\begin{equation}
\labell{eq:FH-K+}
\Z_2\cong\HF_{n+1}^{(Tc+\delta,\,Tc+\eps)}\big(H_+^{(T)}\big)\to
\HF_n^{(Tc-\delta,\,Tc+\delta)}\big(H_+^{(T)}\big)\cong\Z_2
\end{equation}
is an isomorphism.
\end{Lemma}

\begin{proof}
Since \eqref{eq:eps-del0} is satisfied, the results of the previous
section apply, and $x$ and $p$ are the only $T$-periodic orbits of
$\tH_+$ of index $n-1$ or $n$ on $n+1$ with action in the range
$(Tc-\delta,\,Tc+\eps)$.  It is clear that $[p]$ is the generator of
$\HF_n^{(Tc-\delta,\,Tc+\delta)}\big(H_+^{(T)}\big)\cong\Z_2$.
Furthermore, $x$ is the only $T$-periodic orbit of index $n+1$ with
action in $(Tc+\delta,\,Tc+\eps)$ and there are no $T$-periodic orbits
of index $n$ with action in this interval. Hence, the homology
$\HF_{n+1}^{(Tc+\delta,\,Tc+\eps)}\big(H_+^{(T)}\big)$, generated by
$[x]$, is either zero or $\Z_2$. (The former is \emph{a priori}
possible, for in fact there exists a $T$-periodic orbit of index $n+2$
with action in $(Tc+\delta,\,Tc+\eps)$.)  To finish the proof of the
lemma, it is sufficient now to show that the connecting map
\eqref{eq:FH-K+} is onto or, equivalently,
\begin{equation}
\labell{eq:FH0}
\HF_n^{(Tc-\delta,\,Tc+\eps)}\big(H_+^{(T)}\big)=0,
\end{equation}
i.e., the number of Floer anti-gradient trajectories for $\tH_+^{(T)}$
from $x$ to $p$ is odd.

Within $B_R$, the Hamiltonian $T\cdot H_+$ coincides with a standard
bump function $F$ whose variation $C-\min F= T(c-a)$ is greater than
$h(B_R)$ by \eqref{eq:T-delta}. Thus, the assumptions of Section
\ref{sec:bump-W} are satisfied, and $\HF_{n}^{(C-\delta,\,
C+\eps)}(F)=0$. Furthermore, $\tH_+^{(T)}$ agrees with $\tF$ on
$B_R$. Due to our choice of $r$, the condition \eqref{eq:r0} holds and
Lemma \ref{lemma:traj} is applicable. Therefore, $\tF$, and
hence $\tH_+^{(T)}$, have an odd number of Floer anti-gradient
trajectories from $x$ to $p$ contained in $B_R$. Moreover, every Floer
anti-gradient trajectory for $\tH_+^{(T)}$ from $x$ to $p$ is
automatically in $B_R$. This is established by arguing exactly
as in the proof of Lemma \ref{lemma:traj} with $F$ replaced by
$\tH_+^{(T)}$ and using again \eqref{eq:r02}. As a
consequence, $\tH_+^{(T)}$ and $F$ have the same Floer anti-gradient
trajectories from $x$ to $p$, and the total number of such
trajectories is odd. This concludes the proof of \eqref{eq:FH0}
and of the lemma.
\end{proof}

\subsection{The function $H_-$}
Recall that the function $H_+$ and the parameter $\eps>0$ were
fixed above, while $T$ and $\delta$ have been variable.
At this point, we also fix a large period $T$ meeting the requirement
\eqref{eq:T-delta} and such that $T\cdot (c-a)>2\pi r^2$.
Then, condition \eqref{eq:eps-del0} is satisfied if $\delta>0$ is
small, and hence Lemma \ref{lemma:FH1} applies.

In this section, we construct a Hamiltonian $H_-\geq 0$, depending on
$T$, such that $H_-\leq H$ and $H_-(p)=c$, and the connecting map
\begin{equation}
\labell{eq:FH-K-}
\Z_2\cong\HF_{n+1}^{(Tc+\delta_T,\,Tc+\eps)}\big(H_-^{(T)}\big)\to
\HF_n^{(Tc-\delta_T,\,Tc+\delta_T)}\big(H_-^{(T)}\big)\cong\Z_2
\end{equation}
is an isomorphism if $\delta_T>0$ is sufficiently small. 

Recall from Section \ref{sec:prop2-outline}
(see also Definition \ref{def:sympl-deg}) that there exist
\begin{itemize}
\item a loop $\eta^t=\eta^t_i$, $t\in S^1$, of Hamiltonian
diffeomorphisms fixing $p$ and
\item a system of canonical coordinates $\xi=\xi^i$ on a neighborhood
$V$ of $p$
\end{itemize}
such that the Hamiltonian $K=K^i$ generating the flow
$(\eta^t)^{-1}\circ\varphi^t_H$  has a strict local maximum at $p$ and
\begin{equation}
\labell{eq:small}
\max_t\|d^2 (K_t)_p\|_{\xi_p} <\frac{\pi}{ T}.
\end{equation}
Moreover, the loop $\eta$ has identity linearization at $p$, i.e.,
$d(\eta^t)_p=I$ for all $t\in S^1$, and is contractible to $\id$ in
the class of loops with identity linearization at $p$. (See (K3) and
Section \ref{sec:prop2-outline}.) Let $\eta_s$ be a homotopy from
$\eta$ to the identity such that $d(\eta^t_s)_p\equiv I$ and let
$G_t^s$ be the one-periodic Hamiltonian generating $\eta^t_s$ and
normalized by $G_t^s(p)\equiv 0$.  The condition $d(\eta^t_s)_p=I$
is equivalent to that $d^2(G_t^s)_p=0$.

As usual, we normalize $K$ by requiring that $K_t(p)\equiv c$
or, equivalently, by $H=G\# K$. Without loss of
generality, we may also assume that $\bar{V}\subset B_{r_-}$, where
$V$ is the domain of the coordinate system $\xi$ and $B_{r_-}$ is the
ball from the construction of $K_+$; see Section \ref{sec:K+}.
Note that this ball is taken with respect to the original metric and is
not related to $\xi$.

Let $F$ be a bump function, ``centered'' at $p$, with respect to the
coordinate system $\xi$. As in Section \ref{sec:bump-W}, we do not
require $F$ to be supported in $V$, but only constant outside $V$.
Thus, $F\equiv\min F$ on $W\ssminus V$. We may assume
that $\min F<a$. It is also clear that $F$ can be chosen so that
\begin{itemize}

\item $F(p)=c=K(p)$ and $F\leq K$ and, by \eqref{eq:small},
\begin{equation}
\labell{eq:small2}
\|d^2 F_p\|_{\xi_p} <\frac{2\pi}{ T}.
\end{equation}
\end{itemize}
Furthermore, utilizing the condition $d^2(G_t^s)_p=0$ and the
flexibility in the choice of $\min F$ (e.g., making $\min F$ large
negative), we can ensure that
\begin{itemize}
\item $F^s := G^s\# F\leq H_+$ for all $s$.
\end{itemize}

Then $F^s$ is a homotopy beginning with
\begin{equation}
\labell{eq:K-KK+}
H_-:=G^0\#F\leq G^0\#K=H\leq H_+
\end{equation}
and ending with $F^1=F$. Throughout the homotopy, $F^s(p)=c$ and
$F^s\leq H_+$.

The variation of $T\cdot F$, equal to $T(c-\min F)$, is much larger
than $T(c-a)\geq  h(B_{r_-}) \geq h(V)$. Hence, as shown in
Section \ref{sec:bump}, and we have the isomorphism
$$
\Z_2\cong
\HF_{n+1}^{(Tc+\delta_T,\, Tc+\eps)}(T\cdot F)\to
\HF_{n}^{(Tc-\delta_T,\, Tc+\delta_T)}(T\cdot F)
\cong\Z_2,
$$
provided that $\delta_T>0$ is sufficiently small.  Finally note that
$A(G^s)=0$ for all $s$, for $G^s_t(p)\equiv 0$. Therefore, the
functions $F^s$ have equal filtered Floer homology for any period
$T$. In particular, the filtered Floer homology of $H_-^{(T)}$ and of
$F^{(T)}$ are the same and the latter is identical to the
(one-periodic) filtered Floer homology of $T\cdot F$. Thus, we obtain
the desired isomorphism \eqref{eq:FH-K-}.

\subsection{The monotone homotopy map}
By construction, $H_+\geq H_-$. A monotone decreasing homotopy from
$H_+$ to $H_-$ induces maps of filtered Floer homology, which, as is
well known, commute with the maps from the long exact sequence. In
particular, combining the monotone homotopy maps with the connecting
maps \eqref{eq:FH-K+} for $H_+$ and \eqref{eq:FH-K-} for $H_-$, we
obtain the following commutative diagram:
\begin{equation}
\labell{eq:diagram}
\begin{CD}
\Z_2\cong\HF_{n+1}^{(Tc+\delta_T,\,Tc+\eps)}\big(H_+^{(T)}\big) @>\cong>>
\HF_n^{(Tc-\delta_T,\,Tc+\delta_T)}\big(H_+^{(T)}\big) \cong \Z_2\\
@VVV @VV\Psi V\\ 
\Z_2\cong\HF_{n+1}^{(Tc+\delta_T,\,Tc+\eps)}\big(H_-^{(T)}\big)@>\cong>>
\HF_n^{(Tc-\delta_T,\,Tc+\delta_T)}\big(H_-^{(T)}\big) \cong \Z_2
\end{CD}
\end{equation}

To prove the theorem, it is sufficient to show that
the right vertical arrow $\Psi$ in the diagram \eqref{eq:diagram}, i.e.,
the homotopy map
$$
\Psi\colon \HF_n^{(Tc-\delta_T,\,Tc+\delta_T)}\big(H_+^{(T)}\big)
\to
\HF_n^{(Tc-\delta_T,\,Tc+\delta_T)}\big(H_-^{(T)}\big),
$$
is an isomorphism. Indeed, the rows of \eqref{eq:diagram} are
isomorphisms, and hence the left vertical arrow is an isomorphism
whenever $\Psi$ is an isomorphism. Since $H_-\leq H\leq H_+$ by
\eqref{eq:K-KK+}, the left vertical arrow factors as
$$
\HF_{n+1}^{(Tc+\delta_T,\,Tc+\eps)}\big(H_+^{(T)}\big)\to
\HF_{n+1}^{(Tc+\delta_T,\,Tc+\eps)}\big(H^{(T)}\big)\to
\HF_{n+1}^{(Tc+\delta_T,\,Tc+\eps)}\big(H_-^{(T)}\big),
$$
and, as a consequence, the middle group is non-zero as desired.

To show that $\Psi$ is an isomorphism, first observe that since
$F^s\leq H_+$ for all $s$ and $F^0=H_-$ and $F^1=F$, the diagram
$$
\xymatrix{
{\HF_n^{(Tc-\delta_T,\,Tc+\delta_T)}\big(H_+^{(T)}\big)}
 \ar[d]_\Psi\ar[rd] &\\
{\HF_n^{(Tc-\delta_T,\,Tc+\delta_T)}\big(H_-^{(T)}\big)} \ar[r]^{\cong} &
{\HF_n^{(Tc-\delta_T,\,Tc+\delta_T)}\big(F^{(T)}\big)}
}
$$
is commutative, where the horizontal isomorphism is induced by the
isospectral homotopy $F^s$ and the remaining two arrows are monotone
homotopy maps. (See Section \ref{sec:homotopy} and, in particular,
\eqref{eq:diag-iso}.) Recall also that $H_+$ and $F$ are autonomous.
It remains to prove that the diagonal arrow, which can be
identified with
$$
\Z_2\cong\HF_n^{(Tc-\delta_T,\,Tc+\delta_T)}(T\cdot H_+)
\to \HF_n^{(Tc-\delta_T,\,Tc+\delta_T)}(T\cdot F)\cong\Z_2,
$$
is an isomorphism.

Consider a $C^2$-small autonomous perturbation $\hH_+\geq F$ of $H_+$
such that $d^2 (\hH_+)_p$ is negative definite and
$\CS(\hH_+)=\CS(H_+)$.  (It is straightforward to construct
$\hH_+$ by modifying $H_+$ on a neighborhood of
$\bar{B}_{r_-}$.) Then $\hH_+(p)\equiv c$ and $\|d^2 (\hH_+)_p\|$ can be
made arbitrarily small, for $d^2 (H_+)_p=0$.  We take $\hH_+$ such
that $T\cdot \hH_+$ is $C^2$-close to $T\cdot H_+$, and, in
particular, $\|d^2 (T\cdot \hH_+)_p\|$ is small.  Essentially by
definition, the filtered Floer homology of $T\cdot \hH_+$ is
isomorphic to the filtered Floer homology of $T\cdot H_+$, and
it suffices to show that
\begin{equation}
\labell{eq:map}
\Z_2\cong \HF_n^{(Tc-\delta_T,\,Tc+\delta_T)}(T\cdot \hH_+)
\to \HF_n^{(Tc-\delta_T,\,Tc+\delta_T)}(T\cdot F)\cong\Z_2
\end{equation}
is an isomorphism for some small $\delta_T>0$ independent of the
choice of $\hH_+$.

Recall that $Tc$ is an isolated action value of $T\cdot F$ and $T\cdot
H_+$, and hence of $T\cdot\hH_+$, for a generic choice of parameters
of these functions. Fix $\delta_T>0$ meeting the requirement
\eqref{eq:eps-del0} and such that $Tc$ is the only action value of
these Hamiltonians in $(Tc-\delta_T,\,Tc+\delta_T)$. Consider the
linear decreasing homotopy $\hH^s=(1-s)\hH_+ + s F$.  Since both of
the Hessians $d^2 (\hH_+)_p$ and $d^2 F_p$ are negative definite, $d^2
(T\cdot \hH_+^s)_p$ is also negative definite. Thus, $p$ is a
non-degenerate critical point of $\hH_+^s$ for all $s\in [0,\, 1]$ with
$\hH_+^s(p)= c$. Furthermore, by \eqref{eq:small2} and since $\|d^2
(T\cdot \hH_+)_p\|$ is small, $\|d^2 (T\cdot \hH_+^s)_p\|<2\pi$. As a
consequence, $p$ is a uniformly isolated one-periodic orbit of $T\cdot\hH_+^s$;
see, e.g., \cite[pp.\ 184--185]{SZ} or Section
\ref{sec:LFH-LMH}.  By Lemma \ref{lemma:non-zero}, the homotopy map
\eqref{eq:map} is non-zero, and hence an isomorphism.

This concludes the proof of Proposition \ref{prop2} and of Theorem
\ref{thm:main}.  A slightly different proof of the proposition,
although based on the same ideas as the present argument and following
the same line of reasoning, can be found in \cite[Section~5]{GG:gaps}.

\end{document}